\documentclass{article}
\usepackage[utf8]{inputenc}
\usepackage[english]{babel}
\usepackage{mathtools,amsmath,amssymb,amsbsy}
\usepackage{enumitem}
\usepackage{fancyhdr}
\usepackage{csquotes}
\usepackage{geometry}
\usepackage{graphicx}
\usepackage{microtype}
\usepackage{graphicx}
\usepackage{subfigure}
\usepackage{booktabs} % for professional tables

\usepackage{fancyhdr}
\usepackage{ntheorem}
\usepackage{times}  %Required
\usepackage{helvet}  %Required
\usepackage{courier}  %Required
\usepackage{url}  %Required
\usepackage{booktabs}       % professional-quality tables
\usepackage{amsfonts}       % blackboard math symbols
\usepackage{nicefrac}       % compact symbols for 1/2, etc.
\usepackage{microtype}      % microtypography
\newtheorem{theorem}{Theorem}[section]
\newtheorem{remark}{Remark}[section]
\newtheorem{proposition}{Proposition}[section]
\newtheorem{definition}{Definition}[section]
\newtheorem{corollary}{Corollary}[section]
\newtheorem{lemma}{Lemma}[section]

\newtheorem{example}[theorem]{Example}

\theoremstyle{empty}

\newtheorem{refproof}{Proof}[section]

\usepackage[
backend=bibtex,
style=alphabetic,
sorting=ynt
]{biblatex}

\title{Duopoly Investment Problems with Minimally           Bounded Adjustment Costs 
  \newline \newline}

\author{David Mguni\footnote{ Quantitative and Applied Spatial Economic Research Laboratory, University College London, Gower Street, London, WC1E 6BT, UK.}\hspace{1.5 mm}\footnote{Centre for Doctoral Training in Financial Computing \& Analytics, University College London, Gower Street, London, WC1E 6BT, UK. davidmguni@hotmail.com} }

\date{}

\begin{document}

\maketitle
% \label{name}
% \ref{name}
% \nonumber
\begin{abstract}
In this paper, we study two-player investment problems with investment costs that are bounded below by some fixed positive constant. We seek a description of optimal investment strategies for a duopoly problem in which two firms invest in advertising projects to abstract market share from the rival firm. We show that the problem can be formulated as a stochastic differential game in which players modify a jump-diffusion process using impulse controls. We prove that the value of the game may be represented as a solution to a double obstacle quasi-variational inequality and derive a PDE characterisation (HJBI equation) of the value of the game. We characterise both the saddle point equilibrium and a Nash equilibrium for the zero-sum and non-zero-sum payoff games.
\end{abstract}

\textbf{Keywords}: \textit{Impulse control, Stochastic Differential Games, Optimal Stopping, \\Jump-diffusion, Dynkin Games, Verification Theorem, Duopoly, Advertising}.

\section{Introduction}

Over the past three decades, a considerable amount of attention has been dedicated towards modelling the duopolistic advertising problem. The problem is one of finding the optimal dynamic investment strategy for a firm that seeks to maximise cumulative profits over some given time horizon. Each firm uses strategic advertising investments to increase market share. This paper is concerned with an investment problem in which two competing firms make strategic investments over time in order to  maximise their cumulative profits. The paper studies a duopoly environment with future uncertainty. In order to accurately model firm behaviour, we incorporate minimally bounded adjustment costs so that the competing firms incur at least some fixed minimum cost for each investment. This leads to a new description of the advertising investment problem in terms of a non zero-sum stochastic differential game involving impulse controls.

At present, advertising investment models describe firms' investment behaviour with continuous controls --- in particular, it is assumed that the competing firms are able to make infinitesimally fine adjustments to their investment positions that incur arbitrarily small costs. In reality however, advertising investment projects have fixed minimal costs which eliminates the possibility of continual investment since such a strategy would result in singular costs and hence, immediate firm ruin. The presence of fixed minimal costs produces adjustment stickiness (rigidities) since firms adjust their investment positions at discrete points over irregular time intervals. Consequently, the set of feasible investment strategies consists of those in which the firm makes a sequence of investments at selected times along the firm's lifetime.

Despite the relevance of minimally bounded adjustment costs, the literature concerning multiplayer strategic environments with bounded costs remains scarce. Indeed with the exception of the formal mathematical treatment of the zero-sum case presented in \cite{cosso2013stochastic}, models of multiplayer strategic interactions with bounded costs remain limited to environments in which one of the players is allowed to modify the system dynamics continuously (e.g. \cite{yong1994zero,tang1993finite,zhang2011stochastic}).

To account for this, we construct an duopoly investment model in which each firm incurs at least some fixed minimal cost for each advertising investment. To this end, we introduce a non-zero-sum stochastic differential game in which both players use impulse controls to modify the system dynamics. Moreover, in contrast to the models of advertising described above, in order to embed into the description future uncertainty and exogenous economic shocks, we construct a game in which the underlying diffusion process which is allowed to have jumps. 

 The game we study is one in which two players modify a jump-diffusion process using impulse controls in order to maximise some given payoff criterion. We give a PDE characterisation of the value for the game for both zero-sum and non-zero-sum games. The motivation stems from addressing the duopoly advertising investment problem. However, the framework studied and subsequent results derived are general and therefore apply to modelling competitive  multiplayer environments with future uncertainty in which players face fixed adjustment costs. 

We show that the solutions (optimal investment strategies) to the problem can be represented as a solution to a double obstacle problem for a stochastic differential game in which players use impulse controls to modify the system dynamics and give a PDE characterisation of the optimal investment strategies. The theoretical results are complimented by worked examples to illustrate the workings of the theory.
\subsection*{Background Material}

Early versions of the advertising oligopoly problem were formulated as single-player optimal control models in which a controller maximises a payoff extracted from a system whose evolution is governed by some deterministic process. Thus, in the early models of the advertising problem, the influence of competing firms and the effect of future uncertainty derived from market fluctuations and exogenous shocks were neglected (see for example the surveys conducted in \cite{jorgensen1982survey}). To augment the model description, more recent models include a larger repertoire of modelling features, firstly, by modelling the problem as a (two-player) differential game framework the influence of a rival firm can be incorporated into the system --- this approach has yielded considerable descriptive success in modelling the strategic interactions between firms. Following that, \cite{prasad2004competitive} (among others) adopt a stochastic differential game approach to model the problem which accounts for  future uncertainty and random market fluctuations, the inclusion of which has further increased modelling accuracy. We refer the reader to \cite{erickson1995differential} for exhaustive discussions on duopolistic advertising models and to \cite{prasad2004competitive} for a stochastic differential game approach. 

Problems that involve strategic modifications of a controlled dynamic system in competitive environments have attracted much attention over recent years both in theoretical and applied settings. In particular, there is a notable amount of literature on models of this kind in which two players use  continuous controls to modify the system dynamics to satisfy some performance criterion. In the deterministic case, it was shown in \cite{elliott1972existence} and \cite{evans1983differential} that the deterministic differential game admits a value and in fact, the value of the game is a unique solution to a HJBI equation in the viscosity sense. Following on from this, in \cite{fleming1989existence}, the corresponding result was proven for the case in which the system dynamics are stochastic. Indeed, building on the successes of the deterministic cases, the study of stochastic differential game theory has produced significant results and has been successfully applied in various settings within finance and economics. Stochastic differential game theory underpins theoretical models used to prescribe optimal portfolio strategies in a Black-Scholes market (see e.g. \cite{mataramvura2008risk,browne2000stochastic}), descriptions of pursuer-invader dynamics (see e.g. \cite{pachter1981stochastic}) and investment games in competitive advertising (see e.g. \cite{jorgensen1982survey,erickson1995differential,prasad2004competitive}) amongst others. \

If in the differential game, associated to the controllers' modifications to the system dynamics is some fixed minimal cost, the appropriate mathematical framework is a differential game in which the controllers use impulse controls to modify the system dynamics. Impulse control problems are stochastic control models in which the cost of control is bounded below by some fixed positive constant which prohibits continuous control, thus the problem is augmented to one of finding both an optimal sequence of times to apply the control policy, in addition to determining optimal control magnitudes.  Impulse control frameworks therefore underpin the description of financial environments with transaction costs and liquidity risks and more generally, applications of optimal control theory in which the system dynamics are modified by a sequence of discrete actions. We refer the reader to \cite{bensoussan1982controle} as a general reference to impulse control theory and to \cite{korn1999some,jeanblanc1993impulse,palczewski2010finite} for articles on applications.  

Despite the fundamental relevance of fixed minimum costs within economics and financial systems, modelling multi-player competitive economic and financial systems has yet to incorporate the use of impulse control theory. Indeed, unlike in the case of continuous controls for which there is a plethora of studies, with the exception of the zero-sum game studied in \cite{cosso2013stochastic}, stochastic games in which the players use impulse controls remain largely uninvestigated. Deterministic versions of this game were first studied by \cite{yong1994zero,tang1993finite} --- in the model presented in \cite{yong1994zero}, impulse controls are restricted to use by one player and the other uses continuous control. Similarly, in \cite{zhang2011stochastic}, stochastic differential games in which one player uses impulse control and the other uses continuous controls were studied. Using a verification argument, the conditions under which the value of the game (with a single impulse controller) is a solution to a HJBI equation is also shown in \cite{zhang2011stochastic}. In \cite{cosso2013stochastic}, Cosso was the first to study a stochastic differential game in which both players use impulse control using viscosity theory. Thus, in \cite{cosso2013stochastic} it is shown that the game admits a value which is a unique viscosity solution to a double obstacle quasi-variational inequality. In \cite{aid2017nonzero} a nonzero-sum formulation of the game was studied from which a verification theorem that characterises the value function of the game was derived. These results however do not cover the general case of jump-diffusions. \\
\textbf{Contribution}\	

This is the first paper to study the duopolistic advertising problem using a differential game approach with impulsive controls. The paper also makes several theoretical contributions to stochastic differential game theory involving impulse controls. We extend the analyses in \cite{cosso2013stochastic} where stochastic differential games in which both controllers use impulse controls were introduced, to consider i) system dynamics in which the uncontrolled state process is allowed to have jumps ii) non-zero-sum payoff structures. \

We prove verification theorems for both the zero-sum case and the non-zero-sum case in which the appropriate equilibrium concept is a Nash equilibrium. We then generalise the zero-sum payoff structure in the game to a non-zero-sum payoff structure wherein we appeal to a Nash equilibrium as the appropriate equilibrium concept.  Owing to the fact that we use a verification theoretic approach, our proofs differ from those in \cite{cosso2013stochastic}. Both non-zero-sum payoff structures and the inclusion of jumps in the underlying system dynamics serve as important modelling features for applications within finance and economics.
\subsection*{Organisation}

The paper is organised as follows: in Section \ref{section_duopoly}, we introduce our dynamic duopoly model; here we elucidate the key features of the problems and show that the underlying structure is a stochastic differential game. In Section \ref{section_game_description}, we give a description of the stochastic differential game in which impulse controls are used to modify the state process and introduce some of the underlying technology. In Section \ref{section_main_results}, we give a statement of the main results of the paper.  In Section \ref{section_verification_zero_sum}, we prove the main results for zero-sum stochastic differential games with impulse controls culminating in a verification theorem which characterises the value and equilibrium controls of the game. In Section \ref{section_verification_non_zero_sum}, we generalise the results of Section \ref{section_verification_zero_sum} to non zero sum games. Here we characterise the Nash equilibrium of the game. In Section \ref{section_examples}, we study examples that demonstrate the workings of the theory. Here, we give a precise statement of the optimal investment strategies for the duopoly problem. We then give some concluding remarks and lastly, the appendix contains some of the technical proofs from sections \ref{section_verification_zero_sum} - \ref{section_examples}.

We start by providing a description of a generalised duopolistic advertising problem.  
\section{A Duopoly Investment Problem: Dynamic Competitive Advertising}\label{section_duopoly}

The problem is a firm advertising problem in which two firms compete for market share in a duopoly market. Each firm seeks to maximise its long-term profit and to this end, performs investment adjustments, the costs of which are bounded from below. We introduce our model by firstly studying a classical advertising investment model then progressively developing the model to include economic features which include minimally bounded investment costs and exogenous shocks. To fix ideas, as our first case (case I), we consider an environment in which both firms make continuous modifications to their investment positions --- this approach reproduces the \textit{Vidale-Wolfe model of advertising}. Such models do not include fixed minimal costs and assume a zero-sum payoff structure. We refer the reader to \cite{erickson1995differential,jorgensen1982survey} for exhaustive discussions on duopoly advertising investment models and to \cite{prasad2004competitive} for a stochastic differential game approach.

We then consider environments in which adjustment costs are minimally bounded and relax the zero-sum payoff structure in which each firm has its own market share process described by a jump-diffusion process. In this setting, we also embed into our description \textit{cross-over effects} from exogenous shocks from each firm's market share process. The resulting description leads to a new model of dynamic competitive advertising which encapsulates some of the key features of the continuous control model but captures market shocks and cross-over effects between firms. Crucially, the model we introduce requires that each investment incurs at least some fixed minimal cost; this feature of the setup significantly alters the firms' problems and subsequent investment behaviour.

We now give a detailed analysis of each case.   The solution to the model in case II is presented in Sec. \ref{section_examples}.

\subsection*{Overview}

Consider two firms who compete for share of a single market. We refer to the firms as Firm 1 and Firm 2. Both firms seek to maximise profits over some (possibly fixed) horizon by investing in advertising activities in order to increase market share and raise revenue from sales.
\subsection{Case I: Duopoly with Continuous Investments (Review)}

Consider a duopoly in which each firm modifies the advertising investment positions continuously. 
%Therefore, let $S_i (t)\equiv S_i (t,\omega )\in \mathbb{R}\times \Omega ,i\in \{1,2\}$ be a stochastic process on\\ $(\Omega,\mathcal{F},(\mathcal{F}_t )_{t\geq t_0 },\mathbb{P}_0 )$ - a filtered probability space which satisfies the usual conditions of right continuity, augmentation by $\mathbb{P}$-negligible sets and carrying a standard $1$-dimensional $\mathcal{F}_t-$measurable Brownian motion $B$. 
At time $t_0\leq t\leq \tau_S$ each Firm $i\in\{1,2\}$ has a revenue stream $S_i (t)=S_i(t,\omega): \mathbb{R}_{>0}\times \Omega\to \mathbb{R} $ which is a stochastic process and where $\tau_S:\Omega\to  \mathbb{R}_{>0}$ is some common random exit time for the firms' problems. At any point, Firm $i$ makes costly investments of size $u_i\in \mathcal{U}_i$ where $\mathcal{U}_i$ denotes the set of admissible investments for Firm $i$ given $i\in \{1,2\}$. Denote by $M\in \mathbb{R}_{>0}$ the potential market size and by $b_i\in ]0,1]$ the response rate to advertising for Firm $i$, then the revenue stream for Firm $i$ is given by the following expression:
\begin{align}
dS^{t_0  , s_i}_i(t )= b_i  u_i  (t ^-  )\left[M- S_i ^{t_0  , s_i} (t)- S_j ^{t_0  , s_j} (t)\right] M ^{-1}dt  - r_i  S_i ^{t_0  , s_i} (t )dt
+\sigma_i dB (t ),\;\label{marketshareeqn1}
\mathbb{P}-{\rm a.s},&
\\ \text{where}\; i,j\in \{1,2\}\; (i\neq j),
\;\forall (t_0,s_i)\in\mathbb{R}_{>0}\times\mathbb{R}_{>0},&\nonumber
\end{align}
where $s_i\equiv S_i(t_0) \in \mathbb{R}_{>0}$ are the initial sales for Firm $i$; $r_i,\sigma_i\in \mathbb{R}$ are constants that represent the rate at which firm $j\neq i$ abstracts market share and the volatility of the sales process for Firm $i$ respectively. The term $B(t)$ is Brownian motion which introduces randomness to the system.

The cumulative profit for each Firm $i$ is denoted by $\Pi_i$. Each firm seeks to maximise its cumulative profit which consists of its revenue due to sales $h:\mathbb{R}\to \mathbb{R}$ minus its running advertising costs $c_i:\mathbb{R}_{>0}\times \mathbb{R}\to \mathbb{R}$ and lastly, a function of the firm's terminal market share $G:\mathbb{R}\to \mathbb{R}$. 

The profit function for each Firm $i$, $\Pi_i$ is given by the following expression:
\begin{align} 
\Pi_i (t_0,s_i;u_i,u_j )=\mathbb{E}\Bigg[\int_{t_0}^{\tau_S}\left(h(S_i^{t_0,s_i,u_i,u_j } (t))-\left[c_i (t,u_i )-c_j (t,u_j )\right]\right)dt&
+G(S_i^{t_0,s_i,u_i,u_j} (\tau_S ))\Bigg]&. 
\label{profitfunctioncontgame}
\end{align}
Since the market is duopolistic, the payoff structure between the two firms is zero-sum so that the following condition holds:
\begin{equation}
\Pi_1 +\Pi_2 =0. \label{zero_sum_payoff_condition_case_1}
\end{equation}
In light of (\ref{zero_sum_payoff_condition_case_1}), we denote by $\Pi(t,s;u_1,u_2 ):=\Pi_1 (t,s;u_1,u_2 )=-\Pi_2 (t,s;u_1,u_2 ),\\ \forall (t,s)\in\mathbb{R}_{>0}\times \mathbb{R}_{>0},\; \forall u_1\in\mathcal{U}_1, \forall u_2\in\mathcal{U}_2$.

We can now write the dynamic (zero-sum) duopoly problem as:\\
Find $\phi$ and $(\hat{u}_1,\hat{u}_2 )\in \mathcal{U}_1\times \mathcal{U}_2$ s.th.
\begin{align}
\phi(t,s )=\sup_{u_1 \in \mathcal{U}_1}\left(\inf_{u_2 \in \mathcal{U}_2}\Pi^{u_1,u_2} (t,s) \right)=\inf_{u_1 \in \mathcal{U}_1}\left(\sup_{u_2 \in \mathcal{U}_2}\Pi^{u_1,u_2} (t,s) \right)=\Pi^{\hat{u}_1,\hat{u}_2} (t,s),\label{duopolyproblemcontgame}&
\\\forall (t,s)\in \mathbb{R}_{>0}\times \mathbb{R}_{>0},&
\end{align}
where we have used the shorthand $\Pi^{(u_1,u_2 )} (t,s )\equiv \Pi(t,s;u_1,u_2 )$ for any $(t,s)\in\mathbb{R}_{>0}\times\mathbb{R}_{>0}$ and for any $(u_1,u_2)\in\mathcal{U}_1\times\mathcal{U}_2$.

We recognise (\ref{duopolyproblemcontgame}) as a zero-sum stochastic differential game and is a general version of the Vidale-Wolfe advertising model (see for example the differential game extension of the Vidale-wolfe model in \cite{deal1979optimizing}). The stochastic differential game is one in which both players modify the state process using continuous controls. Models of this kind have been used to analyse the strategic interaction within advertising duopoly using a game-theoretic framework \cite{prasad2004competitive}. Using this framework, the behaviour of the firms in the advertising problem can be characterised by computing the optimal policies within a stochastic differential game. 

A feature of this model is that firms are permitted to make infinitely fast investments over the horizon of the problem. This follows since the investment adjustments of each firm are described using continuous controls which allows the firms to make arbitrarily small adjustments to their investment positions (which can incur arbitrarily small costs). Additionally, the zero-sum payoff structure produces a notional transfer of wealth from one firm to the other whenever an advertising investment is made.  

We now present the main duopoly model of the paper:
\subsection*{Case II: Non-Zero-Sum Payoff with Impulse Controls with Jumps}

We now seek a description of the duopoly setting that does not impose a zero-sum payoff structure and removes the ability of firms to perform investment adjustments with arbitrarily small costs. Additionally, we seek a description that accounts for the effect that both firms have on the market which we assume undergoes exogenous shocks. In order to accurately model duopoly investment settings of this kind, it is necessary to embed into the model, fixed minimal bounds to the investment costs which naturally preclude the execution of continuous investment strategies.

To this end, denote by $c:[0,T]\times \mathcal{Z}\to \mathbb{R}$ and by $\chi :[0,T]\times \mathcal{Z}\to \mathbb{R}$ the cost function associated to the advertising investments of Firm $1$ and Firm $2$ respectively where $\mathcal{Z}\subset  \mathbb{R}$ is a given set. Since the firms' investments now incur some minimal cost for each investment, there exist constants $\lambda_1\in \mathbb{R}_{>0}$ and $\lambda_2\in \mathbb{R}_{>0}$ s.th. $c(\tau,z)\geq \lambda_1$ and $\chi (\tau,z)\geq \lambda_2,\;\forall (\tau,z)\in \mathcal{T}\times\mathcal{Z}$. 
%Hence, we now suppose that $(\Omega,\mathcal{F},(\mathcal{F}_t )_{t\geq t_0 },\mathbb{P}_0 )$ - the filtered probability space which satisfies the conditions as case I and now caries both a standard $1-$dimensional $\mathcal{F}_t-$measurable Brownian motion $B$ and compensated Poisson random measure $\tilde{N}$ where $N$ and $B$ are independent. Hence, let $\nu (\cdot)=\mathbb{E}[N(]0,1],\cdot)]$ by a L\'{e}vy measure and where $\tilde{N}(dt,dz)$ is a compensated Poisson measure with $\tilde{N}(dt,dz)\equiv N(dt,dz)-\nu (dz)dt$.\\
  In this case,  continuous investment would result in immediate bankruptcy, each firm must now modify its advertising investment position by performing a discrete sequence of investments over the horizon of the problem. 
  
  Denote by $\mathcal{U}$, the set of admissible investments for Firm $1$, then the sequence of investments $\{\xi_k \}_{k\in \mathbb{N}}$ for Firm $1$  is performed over a sequence of times $\{\tau_k \}_{k\in \mathbb{N}}$. The investment strategy for Firm $1$ is therefore given by a double sequence $u=[\tau_j,\xi_j ]_{j\in \mathbb{N}}\equiv (\tau_1,\tau_2,\ldots;\xi_1,\xi_2,\ldots)\in \mathcal{U}$. Analogously, Firm $2$ has an investment strategy which is the following double sequence: $v=[\rho_m,\eta_m ]_{m\in \mathbb{N}}\equiv (\rho_1,\rho_2,\ldots;\eta_1,\eta_2,\dots)\in \mathcal{V}$, so that Firm $2$ makes a sequence of investments $\{\eta_m\}_{m\in \mathbb{N}}$ which are made over a sequence of times $\{\rho_m \}_{m\in \mathbb{N}}$ where $\mathcal{V}$ is the set of admissible investments for Firm $2$ and where $\xi_1,\xi_2,\ldots,;\eta_1,\eta_2,\ldots,\in \mathcal{Z}$. 
  
  In order to increase its revenue stream, each firm may use advertising investments ($\{\xi_j\}_{j\in\mathbb{N}}$ for Firm 1, $\{\eta_m\}_{m\in\mathbb{N}}$ for Firm 2) to abstract market share which reduces the rival firm's revenue stream. However, increases in either firm's market size expands the economy; this also leads to a higher terminal valuation for both firms which is proportional to the square of the terminal cost (this term is often included in models of duopoly with finite horizon --- see for example \cite{oksendal1998stochastic}). 
  
  The market share processes $S_i$ for each Firm $i$ where $i\in\{1,2\}$ evolve according to the following expressions:
\begin{align}
  &\begin{aligned}
  S_1^{t_0,s_1,u,v}(t)=s_1&+\hspace{-1 mm}\int_{t_0}^{t}\mu_1 S_1^{t_0,s_1,u,v} (r)dr+\sum_{j\geq 1}\xi_j  \cdot 1_{\{\tau_j\leq T \}}  (t)+\hspace{-1 mm}\int_{t_0}^{t}\sigma_{11} (S_1^{t_0,s_1,u,v} (r)) dB_1 (r)
  \\&\begin{aligned}+\int_\mathbb{R}\int_{t_0}^{t}\theta_{11} (S_1^{t_0,s_1,u,v} (r-),z)z\tilde{N}_{11} (dr,dz)+\int_{t_0}^{t}\sigma_{12} (S_1^{t_0,s_1,u,v} (r)) dB_2 (r)&
  \\+\int_\mathbb{R}\int_{t_0}^{t}\theta_{12} (S_1^{t_0,s_1,u,v} (r-),z)z\tilde{N}_{12} (dr,dz).& 	  
   \\
 S_1^{t_0,s_1,u,v} (t_0 )=s_1\in \mathbb{R}_{>0}&,\label{case3stateprocessFirm1}
 \end{aligned}
  \end{aligned}
 \\&
\begin{aligned} 
\hspace{-1.5 mm}S_2^{t_0,s_2,u,v} (t)=s_2&+\hspace{-1 mm}\int_{t_0}^{t}\mu_2 S_2^{t_0,s_2,u,v} (r)dr+\hspace{-1 mm}\sum_{m\geq 1}\eta_m  \cdot 1_{\{\rho_m\leq T \}}  (t)+\int_{t_0}^{t}\hspace{-1 mm}\sigma_{22} (S_2^{t_0,s_2,u,v} (r)) dB_2 (r)
\\&\begin{aligned}
+\int_\mathbb{R}\int_{t_0}^{t}\theta_{22} (S_2^{t_0,s_2,u,v} (r-),z)z\tilde{N}_{22} (dr,dz)
+\int_{t_0}^{t}\sigma_{21} (S_2^{t_0,s_2,u,v} (r)) dB_1 (r)&
\\+\int_\mathbb{R}\int_{t_0}^{t}\theta_{21} (S_2^{t_0,s_2,u,v} (r-),z)z\tilde{N}_{21} (dr,dz)&,	
\\ 
S_2^{t_0,s_2,u,v} (t_0 )=s_2\in \mathbb{R}_{>0}&,
\end{aligned}
\end{aligned}
\label{case3stateprocessFirm2} 
\end{align}
where $\mu_i\in\mathbb{R}$ are given constants, $\tilde{N}_i(ds,dz)\equiv N_i(ds,dz)-\nu(dz)ds$ are compensated random measures\footnote{Recall also that $\nu(\cdot):= \mathbb{E}[N(]0,1],V)]$ for $V\subset\mathbb{R}\backslash\{0\}$.} and $B_i$ are Wiener processes, $i\in\{1,2\}$ and, as before $s_i\equiv S_i^{t_0,s_i,u,v}(t_0) \in \mathbb{R}_{>0}$ are the initial sales for Firm $i$ and $t_0\in[0,T[$ is the start point of the problem. The functions $\sigma_{ii}:\mathbb{R}\to\mathbb{R}$ and $\theta_{ii}:\mathbb{R}\to\mathbb{R}$ represent the \textit{internal} volatility and jump-amplitude for the sales process for Firm $i$ (resp.). The terms $\sigma_{ij}:\mathbb{R}\to\mathbb{R}$ and $\theta_{ij}:\mathbb{R}\to\mathbb{R}$, $i,j\in\{1,2\},i\neq j$, each represent the volatility and jump-amplitudes from the rival firm's sales activities on the firm's own sales process.

Unlike the competitive advertising models that appeal to differential games in which the players' modifications of their investment positions are modelled using continuous controls, 
% (in which firms may make arbitrarily small advertising investments), 
the model now requires that each advertising investment requires at least some fixed minimal cost. 
% Hence here, the firms undertake marketing and advertising projects at discrete points in time in such a way that maximises their cumulative profit.

Our next modification is to relax the zero-sum payoff structure (\ref{zero_sum_payoff_condition_case_1}); we therefore decouple the payoff criterion (\ref{profitfunctioncontgame}) into two profit functions $\Pi_i$ for each Firm $i\in \{1,2\}$. Hence, now Firm $i$ seeks to maximise its running profits over the problem horizon plus a valuation of its terminal market sales. The firms objectives are given by the following expressions:
\begin{align}
&\begin{aligned}
\Pi_1 (t_0,s_1,s_2;u,v)=&\mathbb{E}^{[s_1,s_2 ]}  \Bigg[\int_{t_0}^{\tau_S}e^{-\epsilon r}\left[\alpha_1 S_1^{t_0,s_1,u,v} (r)-\beta_1 S_2^{t_0,s_2,u,v} (r)\right]dr
\\&\hspace{-1 mm}-\sum_{j\geq 1}c_1 (\tau_j,\xi_j ) \cdot 1_{\{\tau_j\leq\tau_S \}} +\gamma_1e^{-\epsilon \tau_S} \left[S_1^{t_0,s_1,u,v } (\tau_S )\right]^2 \left[S_2^{t_0,s_2,u,v} (\tau_S )\right]^2 \Bigg]\label{case3payoffFirm1} 
\end{aligned} 
\\&
\begin{aligned}
\Pi_2 (t_0,s_1,s_2;u,v)=&\mathbb{E}^{[s_1,s_2 ]}  \Bigg[\int_{t_0}^{\tau_S}e^{-\epsilon r}\left[\alpha_2S_2^{t_0,s_2,u,v} (r)-\beta_2 S_1^{t_0,s_1,u,v}(r)\right]dr
\\&\hspace{-6 mm}-\sum_{m\geq 1}c_2 (\rho_m,\eta_m ) \cdot 1_{\{\rho_m\leq\tau_S \}} 
+\gamma_2e^{-\epsilon \tau_S} \left[S_1^{t_0,s_1,u,v} (\tau_S )\right]^2 \left[S_2^{t_0,s_2,u,v} (\tau_S )\right]^2 \Bigg]
\label{case3payoffFirm2}\end{aligned}
\end{align}
where $\tau_S:\Omega \to [0,T]$ is some random exit time at which point the problem ends, $\alpha_i,\beta_i,\gamma_i\in\mathbb{R}$ ($i\in\{1,2\}$) are fixed constants and lastly the constant $\epsilon\in\mathbb{R}_{>0}$ is a common discount factor.

The above model has the following interpretation: the market share $S_i$ of Firm $i$ determines the size of the revenue $\alpha_iS_i$ generated from sales; the parameters $\alpha_i\in\mathbb{R}$ and $\beta_i\in\mathbb{R}$ represent the Firm $i$ revenue per unit sale and the sensitivity of Firm $i$'s activities on Firm $j$'s market share respectively. The terminal quadratic terms capture the fact that increases in either firm's market size \textit{heats up the economy} (expands market opportunities) leading to a higher terminal valuation for \textit{both} firms.

The dynamic duopoly problem is now to characterise the optimal investment policies $(\hat{u},\hat{v})\in\mathcal{U}\times\mathcal{V}$ and to find the functions $\phi_1,\phi_2$ s.th.
\begin{align}
\phi_1(t,s_1,s_2)&=\sup_{u \in \mathcal{U}}\Pi_1 (t,s_1,s_2;u,\hat{v})=\Pi_1 (t,s_1,s_2;\hat{u},\hat{v}),  \label{case3problemFirm1}\\
\phi_2(t,s_1,s_2)&=\sup_{v \in \mathcal{V}}\Pi_2 (t,s_1,s_2;\hat{u},v)=\Pi_2 (t,s_1,s_2;\hat{u},\hat{v}),\qquad	\label{case3problemFirm2}\forall  (t,s_1,s_2)\in \mathbb{R}_{>0}\times \mathbb{R}_{>0}\times \mathbb{R}.
\end{align}

The stochastic differential game therefore involves the use of impulse controls exercised by both players who modify system dynamics to maximise some given payoff function.

In sections \ref{section_verification_zero_sum} and \ref{section_verification_non_zero_sum}, we provide a complete characterisation of the solution to the zero-sum cases and non-zero-sum cases in terms of classical solution to a PDE (respectively). In Section \ref{section_examples} we apply the results of Section \ref{section_verification_non_zero_sum} to characterise the firm's policies for case II. A central component of the proof of the verification theorem is the analysis the players' non-intervention regions. In the zero-sum case,  the opponent's actions produce two changes in the value function: firstly, each impulse action performed by the opponent produces an immediate shift in the value of the state process. This in turn causes indirect changes to the value function since the state process enters as one of its inputs. Secondly, at each intervention, the opponent incurs an intervention cost which, in the zero-sum case represents a transfer of wealth from the opponent to the player --- this produces direct instantaneous changes to the value function. This necessitates reformulating the impulse control problem as a singular control system which has minimally bounded adjustment costs (Lemma \ref{Lemma 3.1.6.}).

We now provide a formal description of the game. We begin by providing a description and the relevant concepts for the zero-sum game which, in Section \ref{section_examples} we shall adapt for the non-zero-sum game:
\section{Description of The Game}\label{section_game_description}

Stochastic differential games are environments in which a number of players interact by altering the dynamics of a stochastic system by strategically selected magnitudes. The players modify the system dynamics in order to maximise some state-dependant payoff criterion over some time horizon, the state process evolves according to some stochastic process that is jointly affected by the players.
\noindent\textbf{Canonical Description}\

The uncontrolled passive state evolves according to a stochastic process $X:[0,T]\times \Omega \to S\subset  \mathbb{R}^p,\; (p\in \mathbb{N})$, which is a jump-diffusion on\\$(\mathcal{C}([0,T]; \mathbb{R}^p ),(\mathcal{F}_{{(0,s)}_{s\in [0,T] }}, \mathcal{F}, \mathbb{P}_{0})$ that is, the state process obeys the following SDE:
\begin{align} 
dX_s^{t_0,x_{0}}=\mu( s,X_s^{t_0,x_{0}} )ds+\sigma(s,X_s^{t_0,x_{0}} )dB_s+\int  \gamma(X_{s^-}^{t_0,x_{0}},z) \tilde{N}(ds,dz)  ,\; X_{t_0}^{t_0,x_{0}}:= x_{0},&
\\\mathbb{P}-{\rm a.s.},\label{uncontrolledstateprocess}
\forall s\in [0,T],\; \forall (t_0,x_{0})\in [0,T]\times S;& 
\end{align}
where $B_s$ is an $m-$dimensional standard Brownian motion, $\tilde{N}(ds,dz)=N(ds,dz)-\nu(dz)ds$ is a compensated Poisson random measure where $N(ds,dz)$ is a jump measure and $\nu(\cdot):= \mathbb{E}[N(1,\cdot)]$ is a L\'{e}vy measure. Both $\tilde{N}$ and $B$ are supported by the filtered probability space and $\mathcal{F}$ is the filtration of the probability space $(\Omega ,\mathbb{P},\mathcal{F}=\{\mathcal{F}_s\}_{s\in [0,T] } )$. We assume that $N$ and $B$ are independent.

We assume that the functions $\mu:[0,T]\times S\to S,\; \sigma:[0,T]\times S\to \mathbb{R} ^{p\times m}$ and $\gamma:\mathbb{R}^p\times \mathbb{R}^l\to \mathbb{R} ^{p\times l}$  are deterministic, measurable functions that are Lipschitz continuous and satisfy a (polynomial) growth condition so as to ensure the existence of (\ref{uncontrolledstateprocess}) \cite{ikeda2014stochastic} (see appendix).

As in \cite{chen2013impulse}, we note that the above specification of the filtration ensures stochastic integration and hence, the controlled jump-diffusion is well defined (this is proven in  \cite{stroock2007multidimensional}).

The generator of $X$ (the uncontrolled process) acting on some function $\phi\in \mathcal{C}^{1,2} (\mathbb{R}^l,\mathbb{R}^p)$ is given by: 
\begin{equation} \mathcal{L}\phi(\cdot,x)=\sum_{i=1}^p   \mu_i (x)    \frac{\partial \phi}{\partial x_i}(\cdot,x)+\frac{1}{2} \sum_{i,j=1}^p  (\sigma \sigma^T )_{ij} (x)    \frac{ \partial^2\phi}{\partial x_i\partial x_{j} }(\cdot,x)+I\phi(\cdot,x)	
\label{generator}
\end{equation}
where $I$ is the integro-differential operator defined by:
\begin{equation} 
I\phi(\cdot,x):= \sum_{j=1}^{l} \int_{\mathbb{R}^p}  \{\phi(\cdot,x+ \gamma^{j}(x,z_j))-\phi(\cdot,x)-\nabla\phi(\cdot,x)  \gamma^{j} (x,z_{j} )\}  \nu_{j} (dz_{j}),\;{\forall  x\in \mathbb{R}^p.}\label{integro_differential_operator} 
\end{equation}

The state process is influenced by a pair of impulse controls  $u=u(s,\omega)\in U;\; s\in [0,T],\omega\in\Omega$ exercised by player I and  $v(s,\omega)\in V;\; s\in [0,T],\omega\in\Omega$ exercised by player II which are each stochastic processes that modify the state process directly. The given sets $U$ and $V$ are convex cones that define the set of admissible controls for player I and player II respectively. 
% The player I control can be written in the form $u=\tilde{f}_1(s,X_s)$ for any $s\in[0,T]$ where $\tilde{f}_1:[0,T]\times S\to \mathfrak{U}$ and $\mathfrak{U}\subset \mathbb{R}^p$ and $\tilde{f}_1$ is some measurable map  w.r.t. $\mathcal{F}$. Analogously, the player II control can be written in the form $v=\tilde{f}_2(s,X_s)$ for any $s\in[0,T]$ where $\tilde{f}_2:[0,T]\times S\to \mathfrak{V}$ and $\tilde{f}_2$ is some measurable map  w.r.t. $\mathcal{F}$  and $\mathfrak{V}\subset \mathbb{R}^p$.  

The player I control is given by $u(s)=\sum_{j\geq 1}\xi_j  \cdot 1_{\{\tau_j\leq T \}}  (s)$ for all $s\in [0,T]$, where $\xi_1,\xi_2,\ldots\in \mathcal{Z}\subset S$ are impulses that are executed at $\mathcal{F}$-measurable stopping times $\{\tau_i\}_{i\in\mathbb{N}}$ where $0\leq t_0\leq \tau_1< \tau_2< \dots <$. Analogously, the player II control is given by $v(s)=\sum_{m\geq 1}\eta_m  \cdot 1_{\{\rho_m\leq T \}}  (s)$  for all $s\in [0,T]$, where $\eta_1,\eta_2,\ldots\in \mathcal{Z}\subset S$ are impulses that are executed at $\mathcal{F}-$measurable stopping times $\{\rho_m\}_{m\in\mathbb{N}}$ where $0\leq t_0\leq \rho_1< \rho_2< \ldots <$. We therefore interpret $\tau_n$ (resp., $\rho_n$) as the $n^{th}$ time at which player I (resp., player II) modifies the system dynamics with an impulse intervention $\xi_n$ (resp., $\eta_n$). We assume that the impulses $\xi_j, \eta_m \in \mathcal{Z}$ are $\mathcal{F}-$measurable for all $m,j \in \mathbb{N}$.  Hence, let us suppose that an impulse $\zeta \in \mathcal{Z}$ determined by some admissible policy $w$ is applied at some $\mathcal{F}-$measurable stopping time $\tau:\Omega \to [0,T]$ when the state is $x'=X^{t_0,x_0,\cdot} (\tau^-)$, then the state immediately jumps from $x'=X^{t_0,x_0,\cdot} (\tau^-)$ to $X^{t_0,x_0,\cdot,w} (\tau)=\Gamma (x',\zeta)$ where $\Gamma :S\times \mathcal{Z}\to S$ is called the impulse response function. 

For notational convenience,  we use $u=[\tau_j,\xi_j ]_{j\geq 1} $ to denote the player I control policy $u=\sum_{j\geq 1}\xi_j  \cdot1_{\{\tau_j\leq T \}}  (s)\in U$.

% The impulses $\xi_j$ (resp., $\eta_m$) are $U-$ valued (resp., $V-$ valued) and are $\mathcal{F}-$measurable for all $j$ (resp., $m$). \

The evolution of the state process with interventions is described by the equation:
\begin{align}
\nonumber X_r^{t_0,x_0,u,v}=x_0&+\int_{t_0}^{r}\mu(s,X^{t_0,x_0,u,v}_s)ds+\int_{t_0}^{r}\sigma(s,X^{t_0,x_0,u,v}_s)dB_s
\\&\nonumber+\sum_{j\geq 1}\xi_j  \cdot 1_{\{\tau_j\leq T \}}  (r)
+\sum_{m\geq 1}\eta_m  \cdot 1_{\{\rho_m\leq T \}}  (r)
\\&+\int_{t_0}^{r}\int\gamma (X^{t_0,x_0,u,v}_{s^-},z) \tilde{N}(ds,dz),\; \mathbb{P}-{\rm a.s},\label{controlledstateprocess}
\\&
\qquad\qquad\quad\qquad\qquad \forall r\in [0,T],\;\forall (t_0,x_0)\in [0,T]\times S,\; \forall u \in U,\forall v \in V.  \nonumber
\end{align}
% where the player I impulse control $u \in U$ and player II impulse control $v \in V $ are given by the following:
% \begin{equation}
% u(s)=\sum_{j\geq 1}\xi_j  \cdot 1_{\{\tau_j\leq \tau_S \}}  (s), \hspace{3 mm} v(s)=\sum_{m\geq 1}\eta_m  \cdot 1_{\{\rho_m\leq \tau_S \}}  (s). 	 
% \label{playerIplayerIIcontrols}
% \end{equation}

Without loss of generality we assume that $X^{t_0,x_0,\cdot}_s=x_0$ for any $s\leq t_0$.

Player I has a gain (or profit) function $J$ which is also a cost function for player II which is given by the following expression:
\begin{align}\nonumber 
J[t_0,x_0;u,v]=\mathbb{E}\Bigg[\int_{t_0}^{\tau_S}\hspace{-1 mm}f (s,X_s^{t_0,x_0,u,v}  ) ds  + \sum_{m\geq 1} c (\tau_m  , \xi_m  )  \cdot 1_{\{\tau_m  \leq  \tau_S  \}}  - \sum_{l\geq 1} \chi(\rho_l  , \eta_l  )  \cdot 1_{\{\rho_l  \leq  \tau_S  \}} 
&\\\nonumber
+G (\tau_S,X_{\tau_S }^{t_0,x_0,u,v})\cdot 1_{\{\tau_S<\infty\}}\Bigg],& 
\\\forall u \in U, \forall v \in V,\;\forall (t_0,x_0)\in [0,T]\times S,& \label{payofffunctionJ}
\end{align}
where $\tau_S:\Omega\to [0,T] $ is some random exit time, i.e. $\tau_S(\omega):=\inf\{s\in [0,T]|X_s^{t_0, x_0  ,\cdot }\in S\backslash A;\;\omega\in\Omega\}$ where $A$ is some measurable subset of $S$, at which point $\tau_S$ the game is terminated. The functions $f:[0,T]\times S\to \mathbb{R} , G:[0,T]\times S\to \mathbb{R}$ are deterministic, measurable, Lipschitz continuous functions which we shall refer to as the \textit{running cost function} and the \textit{bequest function} respectively. The functions $c:[0,T]\times\mathcal{Z}\to\mathbb{R}$ and $\chi:[0,T]\times\mathcal{Z}\to\mathbb{R}$ are deterministic, measurable functions that are the player I and player II intervention cost functions respectively. The assumptions ensure the regularity of the value function (see for example \cite{cosso2013stochastic} and for the single-player case, see for example \cite{davis2010impulse}). We assume that the function $G$ satisfies the condition $\underset{s\to\infty}{\lim}G(s,x)=0$ for any $x\in S$.

In the zero-sum game, the payoff function  $J$ is also the player II cost function which player II minimises.

Note that when either ${U}$ or ${V}$ is a singleton, the game collapses into a classical stochastic impulse control problem with only one player with a value function and solution as that in \cite{oksendal2005applied}.
\subsubsection*{Standing Assumptions}

The results of the paper are built under the assumptions that the functions $f,g,\mu,\sigma,\gamma$ are deterministic, measurable and are Lipschitz continuous and the functions $\mu,\sigma,\gamma$ also satisfy a (polynomial) growth condition (see appendix).

We also make the following assumptions on the cost functions $c,\chi$:\\
A.3.

We assume that there exists constants  $\lambda_c>0$ and $\lambda_\chi >0$  s.th. the following conditions hold $\inf_{\xi \in \mathcal{Z}}c(s,\xi )\geq \lambda_c,$ $\forall$ and $\inf_{\xi \in \mathcal{Z}}\chi (\cdot,\xi )\geq \lambda_\chi $ where $\xi \in \mathcal{Z}$ is a measurable impulse intervention. \\ A.4.

We assume that the cost functions $\chi$  and $c$ are quasi-linear in the impulse inputs. That is to say for all $\mathcal{F}_\tau-$measurable stopping times $\tau\in \mathcal{T}$ and for all $\mathcal{F}_\tau-$measurable impulse interventions $z\in \mathcal{Z}$ we assume the functions $\chi$  and $c$ take the following form:
\begin{equation}
\chi (\tau,z)\equiv a_2  (\tau) z+\kappa_2 \text{ and }c(\tau,z)\equiv a_1(\tau) z+\kappa_1 \label{quasilinearitycosts}\end{equation}
for some constants $\kappa_1,\kappa_2>0$ and some  functions $a_i:[0,T]\to \mathbb{R}$, $i\in\{1,2\}$.

We set $a_2  (\tau)\equiv\lambda_2\in\mathbb{R}$ and $a_1  (\tau)\equiv\lambda_1$ for some $\lambda_1,\lambda_2\in\mathbb{R}^+$ and refer to $\lambda_i$ and $\kappa_i$ as the player $i$ \textit{proportional intervention cost} and \textit{fixed intervention cost} parts (resp.) for $i\in\{1,2\}$.

 Assumption A.3 is integral to the definition of the impulse control problem. Assumption A.4 is necessary to derive an equivalent singular control representation of the problem.

The following definitions are useful:

\begin{definition}\label{Definition 1.2.3.}
Denote by $\mathcal{T}_{(t,\tau')}$  the set of all $\mathcal{F}-$measurable stopping times in the interval $[t,\tau']$, where $\tau'$ is some stopping time s.th. $\tau' \leq T$. If $\tau'=T$ then we denote by $\mathcal{T}\equiv \mathcal{T}_{(0,T)}$. Let $u=[\tau_j,\xi_j]_{j\in\mathbb{N}}$ be a control policy where $\{\tau_j\}_{j\in\mathbb{N}}$ and $\{\xi_j\}_{j\in\mathbb{N}}$ are $\mathcal{F}_{\tau_j}-$ measurable stopping times and interventions respectively, then we denote by $ \mu_{[t,\tau]} (u)$ the number of impulses the controller executes within the interval $[t,\tau]$ under the control policy $u$ for some $\tau \in \mathcal{T}$. 
\end{definition} 

\begin{definition}\label{Definition 1.2.4.}
Let $u\in U$ ($v\in V$) be a player I control (player II control). The impulse control $u$ ($v$) is \textit{admissible} on $[0,T]$ for some $T<\infty$ if either the number of impulse interventions is finite on average that is to say we have that:
\begin{equation}\mathbb{E}[\mu_{[0,T]} (u)]<\infty \qquad (\mathbb{E}[\mu_{[0,T]} (v)]<\infty)	 \end{equation}
or if $\mu_{[0,T]} (u)=\infty\implies$ $\lim_{j\to \infty }\tau_j=\infty$ ($\mu_{[0,T]} (v)=\infty\implies\lim_{m\to \infty }\rho_m=\infty$).
\end{definition}

We shall hereon use the symbol $\mathcal{U}$ (resp., $\mathcal{V}$) to denote the set of admissible controls for player I (resp., player II). For controls $u\in \mathcal{U}$ and $u'\in \mathcal{U}$, we interpret the notion $u\equiv u'$ on $[0,T]$ iff $\mathbb{P}(u=u'$ a.e. on $[0,T])=1$. For the player II controls $v\in\mathcal{V}$ and $v'\in\mathcal{V}$, we interpret the notion $v\equiv v'$ on $[0,T]$ analogously.

\begin{definition}\label{Definition 1.2.5.}
Let $u(s)=\sum_{j\geq 1}\xi_j  \cdot 1_{\{\tau_j\leq T \}}  (s) \in \mathcal{U}$ be a player I impulse control defined over $[0,T]$, further suppose that $\tau:\Omega \to [0,T]$ and $\tau':\Omega \to [0,T]$ are two $\mathcal{F}-$measurable stopping times with $\tau\geq s>\tau'$, then we define the restriction $u_{[\tau',\tau ]} \in \mathcal{U}$ of the impulse control $u(s)$ to be $u(s)=\sum_{j\geq 1}\xi_{\mu_{]t_0,\tau[}(u)+j}  \cdot 1_{\{\tau_{\mu_{[t_0,\tau]}(u)+j} \geq s\geq \tau' \}}  (s)$. 

We define the restriction for the player II control $v(s)=\sum_{m\geq 1}\eta_m  \cdot 1_{\{\rho_m\leq T \}}  (s)$; $v \in \mathcal{V}$  defined over $[0,T]$ analogously so that we define the restriction $v_{[\rho,\rho' ]} \in \mathcal{V}$ of the impulse control $v(s)$ to be $v(s)=\sum_{m\geq 1}\eta_{\nu_{]t_0,\rho[}(u)+m}  \cdot 1_{\{\rho_{\nu_{]t_0,\rho[}(u)+m} \geq s\geq \rho' \}}$ where $\rho\in\mathcal{T}$ and $\rho'\in \mathcal{T}$ are two $\mathcal{F}-$measurable stopping times s.th. $\rho\geq s>\rho'$.
\end{definition}

A player strategy is a map from the other player's set of controls to the player's own set of controls. An important feature of the players' strategies is that they are non-anticipative --- neither player may guess in advance, the future behaviour of other players given their current information. In general, in two player games the player who performs an action first employs the use of a strategy. In particular, the use of strategies affords the acting player the ability to increase its rewards since their action is now a function of the other player's latter decisions.  

Markov controls are those in which the player uses only information about the current state and duration of the game rather than explicitly incorporating information about the other player's decisions or utilising information on the history of the game.  This is in contrast to performing a fixed action or using controls of a Markovian type described above.

In light of the above remarks, limiting the analysis to Markov controls may incur too strong of a restriction on the abilities of the players to perform optimally. It is however, well known that under mild conditions, using Markov controls gives as good performance as an arbitrary $\mathcal{F}$-adapted control (see for example Theorem 11.2.3 in \cite{oksendalapplied2007}). Consequently, in the following analysis we restrict ourselves to Markov controls and hence the player I control takes the form $u=\tilde{f}_1(s,X_s)$ for any $s\in[0,T]$ where $\tilde{f}_1:[0,T]\times S\to \mathfrak{U}$ and $\mathfrak{U}\subset \mathbb{R}^p$ and $\tilde{f}_1$ is some measurable map  w.r.t. $\mathcal{F}$. Analogously, the player II control can expressed in the form $v=\tilde{f}_2(s,X_s)$ for any $s\in[0,T]$ where $\tilde{f}_2:[0,T]\times S\to \mathfrak{V}$ and $\tilde{f}_2$ is some measurable map  w.r.t. $\mathcal{F}$  and $\mathfrak{V}\subset \mathbb{R}^p$.

\begin{definition}
Denote the space of measurable functions by $\mathcal{H}$ and let $\tau\in\mathcal{T}$, we define the [non-local] Player I intervention operator $\mathcal{M}_1:\mathcal{H}\to \mathcal{H}$ acting at a state $X(\tau)$ by the following expression:
\begin{equation}
\mathcal{M}_1 \phi(\tau,X(\tau)):=\inf_{z\in \mathcal{Z}}\left[\phi(\tau,\Gamma (X(\tau^-),z))+c(\tau, z)\cdot 1_{\{\tau\leq T \}}  \right], 	 \label{player1intervention}
\end{equation}
where $\Gamma : S\times \mathcal{Z}\to  S$ is the impulse response function defined earlier.
 We analogously define the [non-local] Player II intervention operator $\mathcal{M}_2:\mathcal{H}\to \mathcal{H}$ at $X(\rho)$ for some $\rho\in\mathcal{T}$ by:
 \begin{equation} \mathcal{M}_2 \phi(\rho,X(\rho)):=\underset{z\in \mathcal{Z}}\sup\hspace{1 mm}\left[\phi(\rho,\Gamma (X(\rho-),z))-\chi (\rho, z) \cdot 1_{\{\rho\leq T\}}  \right].\label{player2intervention}	
 \end{equation}
\end{definition}

Given the remarks of Section \ref{section_game_description}, we now define the value functions of the game:

\begin{definition}\label{Definition 3.1.1.}
The two value functions associated to the game are given by the following expressions:
\begin{align}
\qquad \qquad \qquad \qquad  V^- (t,x)&=\inf_{\hat{u}\in \mathcal{U}}  \sup_{\hat{v}\in \mathcal{V}}  J[t,x;\hat{u},\hat{v} ],\nonumber
\\  V^+(t,x)&=\sup_{\hat{v}\in \mathcal{V}}   \inf_{\hat{u}\in \mathcal{U}}  J[t,x;\hat{u},\hat{v} ],
\\&\nonumber\qquad\qquad
\qquad\qquad\qquad\qquad\qquad\qquad \forall (t,x)\in [0,T]\times S,
\end{align}
where we refer to $V^-$ and $V^+$ as the upper and lower value functions respectively.
\end{definition}

We say that the value of the game exists if $\forall  (t,x)\in [0,T]\times S$ we can commute the supremum and infimum operators in definition 3.7 wherein we can deduce the existence of a function $V\in  \mathcal{C}^{1,2} ([0,T];\mathbb{R}^p )$ with $V(t,x)\equiv V^- (t,x)=V^+ (t,x)$. \

\begin{remark}\label{Remark 3.8.} Suppose the value of the game exists then according to the law of supply and demand, the quantity $\partial_x V(\cdot,x)$ represents the revenue associated with a marginal change in the components thus, $\partial_x V(\cdot,x)$ is the market price when the quantity available is $x$.
\end{remark}
\begin{remark}\label{Remark 3.9.} Suppose the value of the game exists. Then the term $\mathcal{M}_i V(\cdot, x)  ,i\in \{1,2\}$ that is, the non-local intervention operator $\mathcal{M}_i$ acting the value function associated to the game represents the value of the player $i$ strategy that consists of performing the best possible intervention at some given time $s\in [0,T]$ when the state is at $x\in S$, then performing optimally thereafter.

We note however, at any given point an immediate intervention may not be optimal, hence the following inequalities hold:
\begin{align}
\mathcal{M}_1 V(s,x)&\geq V(s,x), 	\label{interventionineqsp1}\\
\mathcal{M}_2 V(s,x)&\leq V(s,x), \quad \forall  (s,x)\in [0,T]\times S. 	\label{interventionineqsp2}    
\end{align} 
\end{remark}
Statements (\ref{interventionineqsp1}) and (\ref{interventionineqsp2}) in fact follow as a direct consequence of the dynamic programming principle. A formal statement and proof for the single impulse controller case is given as Lemma 3.5 in \cite{zhang2011stochastic} for which both proof and results are entirely applicable in the current settings.

Throughout the script we adopt the following standard notation (e.g. \cite{tang1993finite,zhang2011stochastic,chen2013impulse}):

\section*{Notation}

Let $\Omega$  be a bounded open set on  $\mathbb{R}^{p+1}$. Then we denote by:
$\bar{\Omega}$  - The closure of the set $\Omega$.\\
$Q(s,x;R)={{(s',x' ) \in \mathbb{R} ^{p+1}:\max |s'-s|^{\frac{1}{2}}  ,|x'-x|  }<R,s'<s}$. \\
$\partial \Omega$  - The parabolic boundary $\Omega$  i.e. the set of points $(s,x) \in \bar{\mathcal{S}}$ s.th. $R>0, Q(s,x;R)\not\subset\bar{\Omega}$.\\
$\mathcal{C}^{\{1,2\}} ([0, T],\Omega )=\{h \in C^{\{1,2\}} (\Omega ): \partial_s h, \partial_{x_i,x_{j} } h \in C(\Omega )\}$, where  $\partial_s$ and  $\partial_{x_{i}, x_{j}}$ denote the temporal differential operator and second spatial differential operator respectively.\\
$\nabla\phi=(\frac{\partial \phi}{\partial x_1 },\ldots,\frac{\partial \phi}{\partial x_p})$ - The gradient operator acting on some function $\phi \in C^1 ([0,T]\times \mathbb{R}^p)$.\\
$|\cdot|$   - The Euclidean norm to which $\langle x,y \rangle$  is the associated scalar product acting between two vectors belonging to some finite dimensional space.

As in \cite{chen2013impulse}, we use the notation $u={[ \tau_j, \xi_j ]}_{j\geq 1}$ to denote the control policy $u=\sum_{j\geq 1}  \xi_{j}  \cdot 1_{\{ \tau_{j}\leq  T \}}  (s) \in U$ which consists of $\mathcal{F}-$measurable stopping times $\{ \tau_j\}_{j \in \mathbb{N}}$ and $\mathcal{F}$-measurable impulse interventions $\{ \xi_{j}\}_{j \in \mathbb{N}}$. Similarly we use the notation $v={[ \rho_m, \eta_m ]}_{m\geq 1}$ to denote the control policy $v=\sum_{m\geq 1}  \eta_{m}  \cdot 1_{\{ \rho_{m}\leq  T \}}  (s) \in V$ which consists of $\mathcal{F}-$measurable stopping times $\{ \rho_m\}_{m \in \mathbb{N}}$ and $\mathcal{F}$-measurable impulse interventions $\{ \eta_{m}\}_{m \in \mathbb{N}}$.

 \section{Main Results} \label{section_main_results}

We prove two key theoretical results for the game that characterise the conditions for a HJBI equation in both zero-sum and non-zero-sum games. 

	We prove a verification theorem (Theorem \ref{Verification_theorem_for_Zero-Sum_Games with Impulse Control}) for stochastic differential games with a jump-diffusion process and in which the players use impulse controls. In doing so, we prove the following statement:\
	
\begin{theorem}\label{Theorem 4.1.} Suppose that the value of the zero-sum game $V$ exists and that $V\in \mathcal{C}^{1,2} ([0,T],S)\cap \mathcal{C}([0,T],\bar{S} )$, then $V$ satisfies the following double obstacle quasi-variational inequality:
\begin{gather}
\begin{cases}
\begin{aligned}
 \max\Big\{\min\Bigg[-\Big(\partial_s V(t,x)+\mathcal{L}V(t,x)+f(t,x)\Big),&V(t,x)-\mathcal{M}_2 V(t,x)\Bigg],\\
 &V(t,x)-\mathcal{M}_1 V(t,x)\Big\}=0 
 \end{aligned} \\
  V(\cdot,y)=G(\cdot,y), 
  \end{cases}
  \\\nonumber\qquad\qquad\qquad\qquad\qquad\qquad\qquad\qquad\qquad\qquad\qquad\qquad\forall y \in S,\; \forall (t,x)\in [0,T]\times S,
 \end{gather}
where $\mathcal{M}_i$ is the [non-local] player $i$ intervention operator for $i\in\{1,2\}$, $T\in\mathbb{R}_{>0}$ is the time horizon of the game and $S\subset \mathbb{R}^p$ is a given set. 

Moreover, denote by $\hat{u}:=[\hat{\tau}_j,\hat{\xi}_j ]_{j\in \mathbb{N}}$ and $\hat{v}:=[\hat{\rho}_m,\hat{\eta}_m ]_{m\in \mathbb{N}}$ the equilibrium controls, then $V$ satisfies the following expression:
\begin{equation} 
 \Delta_{\hat{\xi}_j} V(\hat{\tau}_j,\cdot)=-c(\hat{\tau}_j,\hat{\xi}_j) ,\label{marginalcostcondition}\end{equation}
where $\Delta_z\phi(\cdot,X(\tau)):=\phi(\cdot,\Gamma (X(\tau^-),z))-\phi(\cdot,X(\tau^-))+\Delta_NX(\tau) $ given some $\mathcal{F}_\tau-$measurable intervention $z\in\mathcal{Z}$ and where $\Delta_N X(\tau)$ denotes a jump at some $\mathcal{F}_\tau-$measurable stopping time $\tau\in\mathcal{T}$ due to $\tilde{N}$.
\end{theorem}

For the non zero-sum payoff case with running cost functions $f_i$ and terminal cost functions $G_i$ for $i\in\{1,2\}$, we have the following result:

\begin{theorem}\label{Theorem 4.2.} 
Denote by $\phi_i$ the player $i$ value function for the non zero-sum game for $i\in \{1,2\}$ then the value functions $\phi_i$ satisfy the following quasi-variational inequalities:
\begin{align}
\begin{cases} \max\{-(\partial_s \phi_i(t,x)+\mathcal{L}\phi_i(t,x)+f_i(t,x)),\phi_i(t,x)-\mathcal{M}_i \phi_i(t,x)\}=0,  \\
  \phi_i(\tau_S,y)=G_i(\tau_S,y),\qquad\forall y \in \mathbb{R},\; (t,x)\in [0,T]\times \mathbb{R},\; i\in\{1,2\}. 
  \end{cases}
  \end{align}
  \end{theorem}

The following results characterise the optimal investment policies for the duopoly problem described in Section \ref{section_duopoly}:

\begin{theorem}\label{Theorem 4.3.} Suppose that the market share $X_i$ of Firm $i$, ($i \in \{1,2\}$) evolves according to (\ref{case3stateprocessFirm1}) - (\ref{case3stateprocessFirm2}) and let the firm payoff functions be given by (\ref{case3payoffFirm1}) - (\ref{case3payoffFirm2}), then the sequence of optimal investments $\hat{u}=[\hat{\tau}_j,\hat{\xi}_j ]_{j\in \mathbb{N}}\equiv\sum_{j\geq 1}\hat{\xi}_j  \cdot 1_{\{\hat{\tau}_j\leq \tau_S \}}  (s)$ for Firm 1 is characterised by the investment times $\{\hat{\tau}_j\}_{j\in\mathbb{N}}$ and investment magnitudes $\{\hat{\xi}_j\}_{j\in\mathbb{N}} $ where $[\hat{\tau}_j,\hat{\xi}_j ]_{j\in \mathbb{N}}$ are constructed via the following expressions:
\renewcommand{\theenumi}{\roman{enumi}}
\begin{enumerate}
\item $\hat{\tau}_0\equiv t_0\text{ and  }\hat{\tau}_{j+1}=\inf\{s>\tau_j;X_1 (s)\leq x^{\star}_1\}\wedge \tau_S, \hspace{1 mm}\forall  v \in \mathcal{V}$.
\item $\hat{\xi}_j=\hat{x}_1-X_1(\hat{\tau}_j)$.
\end{enumerate}

Similarly for Firm 2, the optimal sequence of investments $\hat{v}:=[\hat{\rho}_m,\hat{\eta}_m ]_{m\in \mathbb{N}}=\sum_{m\geq 1}\hat{\eta}_m  \cdot 1_{\{\hat{\rho}_m\leq \tau_S \}}  (s)$ is given by:
\renewcommand{\theenumi}{\roman{enumi}}
\begin{enumerate}
\item$\hat{\rho}_0\equiv t_0\text{ and  }\hat{\rho}_{m+1}=\inf\{s>\rho_m;X_2 (s)\leq x^{\star}_2 \}\wedge \tau_S$,  $  \hspace{1 mm}\forall  u \in \mathcal{U}$,
\item $\hat{\eta}_m=\hat{x}_2-X_2(\hat{\rho}_m)$, 
\end{enumerate}
where the quadruplet $(x_1^{\star},x_2^{\star},\hat{x}_1,\hat{x}_2)$ is determined by the following equations ($i\in\{1,2\}$):
\begin{align}
C_1r_{1,i}e^{r_{1,i}x_i^{\star}}&+C_2r_{2,i}e^{r_{2,i}x_i^{\star}}+\frac{\alpha_i}{\epsilon}=\lambda_i\label{Riccconditions1.0}\\
C_1r_{1,i}e^{r_{1,i}\hat{x}_i}&+C_2r_{2,i}e^{r_{2,i}\hat{x}_i}+\frac{\alpha_i}{\epsilon}=\lambda_i\label{Riccconditions2}\\
C_1(e^{r_{1,i}x_i^{\star}}-e^{r_{1,i}\hat{x}_i})&+C_2(e^{r_{2,i}x_i^{\star}}-e^{r_{2,i}\hat{x}_i})=-\kappa_i +\Big(\lambda_i-\frac{\alpha_i}{\epsilon}\Big)(x^{\star}_i-\hat{x}_i),
\label{Riccconditions3.0}\end{align}
where $C_1$ and $C_2$ are endogenous constants whose values are determined by (\ref{Riccconditions1.0}) - (\ref{Riccconditions3.0}), $\lambda_i$ and $\kappa_i$ are the Firm $i$ proportional and fixed intervention costs (respectively),  $\alpha_i$ is the Firm $i$ margin parameter, $\epsilon$ is the discount rate and the values $r_{1,i}$ and $r_{2,i}$ are roots of the equation:
\begin{equation}
q(r_{k,i}):=\frac{1}{2}\sigma_{ii}^2r_{k,i}^2+\mu_i r_{k,i}-\epsilon+\int_{\mathbb{R} }\{e^{r_{k,i}\theta_{ij}z}-1-\theta_{ij}r_{k,i}z\}\nu_j(dz),\label{rrootequation}
\end{equation}
for $i,j,k\in\{1,2\}$.
\end{theorem}

Theorem \ref{Theorem 4.3.} says that each firm performs a sequence of investments over the time horizon of the problem. The decision to invest is determined by the firm's market share process --- in particular, at the point at which Firm $i$'s share of the market falls below the level $x^{\star}_i$, then the firm performs an investment in order to raise its market share to $\hat{x}_i$, where the fixed values $\hat{x}_i$ and $x^{\star}_i$ are determined by the given parameters $\lambda_i,\lambda_j,\kappa_i,\kappa_j,\alpha_i,\epsilon$ via (\ref{Riccconditions1.0}) - (\ref{Riccconditions3.0}). In particular, each Firm $i$ seeks to retain a market share of at least $x^{\star}_i$, where $x^{\star}_i$ is a quantity determined by the size and influence of both firms. 
% At any point Firm $i$'s market share falls below $x^{\star}_i$ the firm immediately reacts by performing an investment in order to raise its market share to $\hat{x}_i$ and in doing so abstracting market share from its competitors. 
Therefore, if $S$ is the total size of the market the value $S-x^{\star}_i$ represents the maximum level of market share that Firm $i$ is prepared to cede to the rival firm.

In summary, each firm observes its own market share and only intervenes at the points at which the firm's market share has fallen below some fixed level. 
% At this point, the firm performs an advertising investment in order to raise the market share to within some prefixed levels. For each firm, both the minimum market share and the investment magnitudes are determined by the firm's size and the responsiveness of the market to advertising investments of both firms.
Each firm's intervention policy is reactant to the investment and subsequent market acquisition of the other firm, each firm therefore reacts by performing the best sequence of response investments to the other firm's investment strategy.

The following corollary follows directly from Theorem \ref{Theorem 4.3.} and establishes when each firm performs investments under the optimal Nash equilibrium strategy:\

\begin{corollary}\label{Corollary 4.4.} The sample space splits into three regions: a region in which Firm $i$ performs an advertising investment --- $I_1$, a region in which Firm 2  performs an advertising investment --- $I_2$ and a region in which no action is taken by either firm $I_3$. Moreover, the three regions are characterised by the following expressions:
\begin{align*}
&I_j=\{x< x_j^{\star}|x,x_j^{\star}\in\mathbb{R}\}, \qquad i,j\in \{1,2\},\\&
I_3=\{x\geq x^{\star}_i\wedge x_j^{\star}|x,x_i^{\star},x_j^{\star}\in\mathbb{R}\},
\end{align*}
where the $x^{\star}_i, x^{\star}_j, i,j \in \{1,2\}$ are determined by (\ref{Riccconditions1.0}) - (\ref{Riccconditions3.0}).
\end{corollary}
The following result characterises the value function for each firm:
\begin{proposition}\label{firm_value_function_ch4}
The value function $\phi_i(t,x_1,x_2):\mathbb{R}_{>0}\times\mathbb{R}_{>0}\times\mathbb{R}_{>0}$ for each Firm $i$ is given by the following expression:
\begin{align}
\phi_i(t,x_1,x_2)=e^{-\epsilon t}\Big\{C_1(e^{r_1x_1}+e^{r_1x_2})+C_2(e^{r_2x_1}+e^{r_2x_2})+\frac{\alpha_i}{\epsilon}x_1&-\frac{\beta_i}{\epsilon}x_2
\\&\begin{aligned}+\frac{1}{2\epsilon^2}(\mu_i\alpha_i-\mu_j\beta_i)\Big\},&
\\i\neq j,\nonumber i,j\in\{1,2\},&
\label{psi_solution_proposition}
\end{aligned}
\end{align}
where $C_1,C_2,r_1,r_2$ are endogenous constants.
\end{proposition}

Proposition \ref{firm_value_function_ch4} provides a full characterisation of each firm's value function which in turn, quantifies each firm's future expected payoff. As we show in the appendix, the endogenous constants can be recovered by approximating the solutions to a system of simultaneous equations. 

% To our knowledge, this is the first time a jump-diffusion process within a stochastic differential game in which the players use impulse controls to modify the state process has been considered.

In the following analysis, we use the results of the stochastic differential game involving impulse controls and a non zero-sum payoff to solve our model of investment duopoly (case II) presented in the paper.

To our knowledge, this paper is the first to deal with a jump-diffusion process within a stochastic differential game in which the players use impulse controls to modify the state process.

\section{A HJBI Equation for Zero-Sum Stochastic Differential Games with Impulse Controls.} \label{section_verification_zero_sum}

In this section, we give a verification theorem for the value of the game therefore giving conditions under which the value of the game is a solution to the HJBI equation. 

To accommodate the influence of impulses exercised by player II on the value function, it is necessary to reformulate the problem as a singular impulse control problem for which we appeal to Lemma \ref{Lemma 3.1.6.}.

\subsection*{Preliminaries}

First, we give some preliminary results which are necessary to prove the main results of the section.

\begin{lemma}\label{Lemma 3.1.7.}
Let $(\tau,x)\in \mathcal{T}\times S$ and let $V\in \mathcal{H}$, then the sets $\Xi_1$ and $\Xi_2$ defined by:
\begin{gather}
\Xi_1 (\tau,x):=\big\{\xi \in \mathcal{Z}:\mathcal{M}_1 V(\tau^-,x)=V(\tau^-,x+\xi )+c(\tau, \xi )\cdot  1_{\{\tau\leq T \}}  \big\},\\
\Xi_2 (\tau,x):=\big\{\xi \in \mathcal{Z}:\mathcal{M}_2 V(\tau^-,x)=V(\tau^-,x+\xi )-\chi (\tau, \xi ) \cdot 1_{\{\tau\leq T \}}  \big\},
\end{gather}
are non-empty.\end{lemma}

The proof of the lemma is essentially that given as the proof of Lemma 3.7 in \cite{chen2013impulse} with little adaptation --- we therefore omit the proof of lemma here.

\begin{definition}[\cite{guo2008solving}]\label{Definition 3.1.5.}
Let $I\subset \mathbb{R}$ be an open (and possibly unbounded) interval and denote its closure by $\bar{I}$. Suppose $x\in  \bar{I}$, then an admissible singular control is a pair $(\nu_s^+,\nu_s^- )_{s\geq t}$ of $\mathcal{F}-$adapted, non-decreasing c\`{a}dl\`{a}g processes s.th. $\nu ^+ (0)=\nu ^- (0)=0,\; X_s^{t_0,x_0,\cdot}:=x_0+\nu_s^+-\nu_s^-$ and $d\nu ^+, d\nu ^-$ are supported on disjoint subsets.
\end{definition}
The following result demonstrates that general impulse control problems can be represented as a singular control problem, in particular it shows that the game (\ref{payofffunctionJ}) can be represented as a game of singular control.
\begin{lemma}\label{Lemma 3.1.6.}
Let $\nu\equiv(\nu_1 (s),\nu_2(s)):[0,T]\times\Omega\times[0,T]\times\Omega\to\mathbb{R}^{2p}$ be a pair of adapted finite variation c\`{a}dl\`{a}g processes with increasing components and let $\Theta_i\in\mathcal{H}, i \in \{1,2\}$ be a pair of functions which satisfy conditions A.2.1 - A.2.2. The impulse control problem with cost functions $c$ and $\chi$ given by $c(\tau_k,\xi_k)\equiv\lambda_1\xi_k+\kappa_1$ and $\chi(\rho_m,\eta_m)\equiv\lambda_2\eta_m+\kappa_2$ for player I and player II respectively where $\xi,\eta\in\mathcal{Z};\; \lambda_i,\kappa_i\in\mathbb{R}_{>0},\; i\in\{1,2\}$, is equivalent to the following singular control problem:

Find $\phi\in\mathcal{H}$ and $\hat{\nu}$ s.th.
    \begin{align} 
    \phi(t_0,x_0)=\inf_{\nu_1}\sup_{\nu_2}J[t_0,x_0;\nu ]=J[t_0,x_0;\hat{\nu} ],
\end{align}
where
    \begin{align}J[t_0,x_0;\nu ]=\mathbb{E}\Bigg[\int_{t_0}^{\tau_S}f (s,X_s^{t_0,x_0,\nu} )ds+\int_{t_0}^{\tau_S}\Theta_1(s)d\nu_1 (s)&-\int_{t_0}^{\tau_S}\Theta_2(s)d\nu_2 (s)
    \\&\begin{aligned}+G(\tau_S,X_{\tau_S}^{t_0,x_0,\nu} )\cdot 1_{\{\tau_S<\infty\}}\Bigg],&\label{singcontrolproblem}
\\ \forall (t_0,x_0)\in [0,T]\times S,& \nonumber
\end{aligned}
\end{align}
and where the state process $X$ evolves according to the following SDE:
\begin{align}
X^{t_0,x_0,\nu}_r=x_0+\int_{t_0}^{r}\mu(s,X^{t_0,x_0,\nu}_s)ds&+\int_{t_0}^{r}\sigma(s,X^{t_0,x_0,\nu}_s)dB_s+\nu (r)
\\&\begin{aligned}+\int_{t_0}^{r}\int\gamma (X^{t_0,x_0,\nu}_{s-},z) \tilde{N}(ds,dz),& \label{singularcontrolledstateprocess}
\\\mathbb{P}-{\rm a.s}., \hspace{1 mm}\forall (t_0,x_0)\in [0,T]\times S.\hspace{-5 mm}&
\end{aligned}
\end{align}
\end{lemma}

We will use Lemma \ref{Lemma 3.1.6.} in order to prove a verification theorem for the stochastic differential game which has a zero-sum payoff structure.

We defer the proof of Lemma \ref{Lemma 3.1.6.} to the appendix.

The following theorem characterises the conditions in which the value of the game satisfies a HJBI equation:

\begin{theorem}[Verification theorem for Zero-Sum Games with Impulse Control]\label{Verification_theorem_for_Zero-Sum_Games with Impulse Control}

% Let $\tau$ be some $\mathcal{F}-$measurable stopping time and denote by $\hat{X}(\tau)=X(\tau^- )+\Delta_N X(\tau)$ where $\Delta_N X(\tau)$ denotes a jump at some $\mathcal{F}_\tau-$measurable time $\tau$ due to $\tilde{N}$. Denote also by $\Delta_z\phi(X(\tau)):=\phi(\Gamma (X(\tau^-),z))-\phi(X(\tau^-))+\Delta_NX(\tau)$ where $\tau\in\mathcal{T}$ and $z\in\mathcal{Z}$ is some $\mathcal{F}_\tau-$measurable stopping time and intervention (resp.).

Suppose that the value of the game $V$ exists and that $V\in \mathcal{C}^{1,2} ([0,T],S)\cap \mathcal{C}([0,T],\bar{S} )$. Suppose also that there exists a function $\phi\in \mathcal{C}^{1,2}([0,T],S)\cap \mathcal{C}([0,T],\bar{S} )$ that satisfies technical conditions (T1) - (T4) and the following conditions $\forall s \in [0,T]$: 
	 \renewcommand{\theenumi}{\roman{enumi}}
 \begin{enumerate}[leftmargin= 8.5 mm]
\item\label{ii_ch4_intervention_ineq} $\phi\leq \mathcal{M}_1 \phi$ in $S$ and $\phi\geq \mathcal{M}_2 \phi$ in $S$ where $D_1$ and $D_2$ are defined by:
$D_1=\{x\in S;\phi(\cdot,x)<\mathcal{M}_1 \phi(\cdot, x)\}$ and $D_2=\{x\in S;\phi(\cdot, x)>\mathcal{M}_2 \phi(\cdot,x)\}$
where we refer to $D_1$ (resp., $D_2)$ as the player I (resp., player II) continuation region.
\item \label{iii_hjbi_p1_ineq} $\frac{\partial \phi}{\partial s}+\mathcal{L}\phi(s,X^{\cdot,\hat{u} ,v  } (\cdot))+f(s,X^{\cdot,\hat{u} ,v  } (\cdot))\leq 0    $,  $\forall  v \in \mathcal{V}$ on $S\backslash{\partial D_2}$.
\item \label{iv_hjbi_p2_ineq} $\frac{\partial \phi}{\partial s}+\mathcal{L}\phi(s,X^{\cdot,u ,\hat{v}} (\cdot))+f(s,X^{\cdot,u ,\hat{v}} (\cdot))\geq 0   $,  $\forall  u \in \mathcal{U}$ on $S\backslash{\partial D_1}$.
\item\label{v_ch4_continuity_equation_equality} $\frac{\partial \phi}{\partial  s}+\mathcal{L}\phi(s,X^{\cdot,\hat{u},\hat{v} } (\cdot))+f(s,X^{\cdot,\hat{u},\hat{v} } (\cdot))=0$ in $D\equiv D_1\cap D_2$ in   $S$.   
\item\label{vi_ch4_terminal_limit} $X^{\cdot,u,v} (\tau_S )\in \partial S$, $\mathbb{P}-$a.s. on $\{\tau_S<\infty \}$ and $\phi(s,X^{\cdot,u,v} (\cdot))\to G(\tau_S,X^{\cdot,u,v} (\tau_S )) \cdot 1_{\{\tau_S<\infty \}}$  as $s\to \tau_S^-\;  \mathbb{P}-$a.s.,$\forall  x\in S, \forall u \in \mathcal{U},\; \forall v \in \mathcal{V}$. 
\item \label{vii_ch4_terminal_limit} $\hat{\xi}_k\in \arg\hspace{-0.55 mm}\inf_{z\in \mathcal{Z}}\{\phi(\tau_k^-,\Gamma (x,z))+c(\tau_k,z)\}$ is a Borel Measurable selection and similarly, $\hat{\eta}_j\in \arg\hspace{-0.55 mm}\sup_{z\in \mathcal{Z}}\{\phi(\rho_j^-,\Gamma (x,z))-\chi (\rho_j,z)\}$ is a Borel Measurable selection $\forall  x\in S;$ and for any $\tau_1,\tau_2,\ldots;\rho_1,\rho_2,\ldots\in\mathcal{T}$.
\end{enumerate}

Put $\hat{\tau}_0\equiv t_0$ and define $\hat{u}:=[\hat{\tau}_j,\hat{\xi}_j ]_{j\in \mathbb{N}}$ inductively by: $\hat{\tau}_{j+1}=\inf\{s>\tau_j;\; X^{\cdot,\hat{u} ,v} (\cdot)\notin D_1 \}\wedge \tau_S,\; \forall  v \in \mathcal{V}$. Similarly, put $\hat{\rho}_0\equiv t_0$ and define  $\hat{v}:=[\hat{\rho}_m,\hat{\eta}_m ]_{m\in \mathbb{N} } $ inductively by $\hat{\rho}_{m+1}=\inf\{s>\rho_m;\; X^{\cdot,u ,\hat{v}} (\cdot)\notin D_2 \}\wedge \tau_S, \; \forall  u \in \mathcal{U}$.\
\renewcommand{\theenumi}{\roman{enumi}}\begin{enumerate}\setcounter{enumi}{6}
\item $\label{verif_1_singular_cond}\qquad\qquad\Delta_z\phi(\hat{\rho}_m),X(\hat{\rho}_m))=\chi(\hat{\rho}_m,z)  \text{ and }  \Delta_{\hat{\xi}_j}\phi(\hat{\tau}_j,X(\hat{\tau}_j))=-c(\hat{\tau}_j,z),\;\forall z\in\mathcal{Z}$.
  \end{enumerate}
Then 
\begin{align}
\phi(t,x)=J[t,x;\hat{u},\hat{v}]=\inf_{u\in \mathcal{U}}  \sup_{v\in \mathcal{V}} J[t,x;u,v]=  \sup_{v\in \mathcal{V}}   \inf_{u\in \mathcal{U}}  J[t,x;u,v ],\qquad \forall (t,x)\in [0,T]\times S.
\end{align}
\end{theorem}
% Theorem \ref{Verification_theorem_for_Zero-Sum_Games with Impulse Control} says that if a sufficiently smooth candidate solution to the HJBI equation in \eqref{v_ch4_continuity_equation_equality} can be found then the value function of the game is characterised in terms of a solution to the PDE in \eqref{v_ch4_continuity_equation_equality}. 
Theorem \ref{Verification_theorem_for_Zero-Sum_Games with Impulse Control} provides a full characterisation of the value for the game. In particular, the theorem states that in equilibrium, both players play QVI controls and, provided that the value of the game exists and is sufficiently smooth to apply the Dynkin formula and, should a solution to the HJBI equation (\ref{v_ch4_continuity_equation_equality}) exist, then the value of the function coincides with the HJBI equation solution.

From Theorem \ref{Verification_theorem_for_Zero-Sum_Games with Impulse Control}, we also see that the sample space splits into three regions that consist of a continuation region, in which neither player performs an intervention and intervention regions for each player within which, the player performs an impulse execution. In particular, we have the following corollary:

\begin{corollary}\label{Corollary 6.2.} The sample space splits into three regions which, when playing their equilibrium strategies, represent a region in which player I executes interventions $I_1$, a region in which player II executes interventions $I_2$, and a region $I_3$ in which no action is taken by either player; moreover the three regions are characterised by the following expressions:
\begin{align*}
&I_1=\{(t,x)\in [0,T]\times S: V(t,x)=\mathcal{M}_1 V(t,x),\;\mathcal{L}V(t,x)+f(t,x)\geq 0\},
\\
&I_2=\{(t,x)\in [0,T]\times S: V(t,x)=\mathcal{M}_2 V(t,x),\mathcal{L}V(t,x)+f(t,x)\geq 0\}  ,
\\ 
&\begin{aligned}
I_3= \{(t,x)\in  [0,T]\times S: V(t,x)<\mathcal{M}_1 (t,x),&V(t,x)>\mathcal{M}_2 V(t,x);
\\&\;\;\;\;\mathcal{L}V(t,x)+f(t,x)=0\}
\end{aligned}
\end{align*}
\end{corollary}

\begin{refproof}[Proof of Theorem \ref{Verification_theorem_for_Zero-Sum_Games with Impulse Control}]

In the following, we make the distinction between the jumps due to the players' impulse controls and the jumps due to $\tilde{N}$. Indeed, for any $\tau\in\mathcal{T}$, we denote by $\hat{X}(\tau):=X(\tau^- )+\Delta_N X(\tau)$ where $\Delta_N X(\tau)$ is the jump at $\tau$ due to $\tilde{N}$ where $\Delta_N X(s)=\int  \gamma(X(s-),z) \tilde{N}(ds,dz)$ and $\tilde{N}(ds,dz)=\tilde{N}(s,dz)-\tilde{N}(s,dz)$. \

Similarly, given an impulse $\xi \in \mathcal{Z}$ (resp., $\eta\in \mathcal{Z}$) exercised by player I (resp., player II), we denote the jump induced by the player I (resp., player II) impulse by $\Delta_\xi$  (resp., $\Delta_\eta$). That is, for any $\tau\in\mathcal{T}$, we define $\Delta_\xi  \phi(\tau,X^{t_0,x_0,u,v} (\tau)):=\phi(\tau^-,\Gamma (X^{t_0,x_0,u,v} (\tau^-),\xi ))-\phi(\tau^-,X^{t_0,x_0,u,v} (\tau^-))+\Delta_N \phi(\tau,X^{t_0,x_0,u,v} (\tau))$ to be the change in $\phi$ due to the player I impulse $\xi \in \mathcal{Z}$ where $\Gamma : S\times Z\to  S$ is the impulse response function. 

We define $\Delta_\eta$ analogously so that $\Delta_\eta \phi(\rho,X^{t_0,x_0,u,v} (\rho))\\:=\phi(\rho^-,\Gamma (X^{t_0,x_0,u,v} (\rho^-),\eta))-\phi(\rho^-,X^{t_0,x_0,u,v} (\rho^-))+\Delta_N \phi(\rho,X^{t_0,x_0,u,v} (\rho))$ is the change in $\phi$ due to the player II impulse $\eta\in \mathcal{Z}$ at some intervention time $\rho\in \mathcal{T}$.

To prove the theorem, we use a singular control representation of the combined impulse controls for each player. For our first case, we define $\nu$ by  $\nu(s)\equiv\eta(s)+\xi(s)$ so that $\nu$ is a process consisting of the combined player I and player II controls. Note that by Lemma \ref{Lemma 3.1.6.}, we have the following equivalences $\forall r \in [0,T]$:
\renewcommand{\theenumi}{\alph{enumi}}
 \begin{enumerate}[label={(\alph*)}]
\item $\xi(r)=\sum_{m=1}^{\mu_{[t_0,r]}}\xi_j\cdot1_{\{\tau_j\leq T\}}$
\item $\eta(r)=\sum_{m=1}^{\mu_{[t_0,r]}}\eta_m\cdot1_{\{\rho_m\leq T\}}
$\item $\int_r^{r'}d\xi(s)=\sum_{j=1}^{\mu_{[r,r']}(v)}c(\tau_j,\xi_j)\cdot1_{\{\tau_j\leq T\}}$
\item $\int_r^{r'}d\eta(s)=\sum_{m=1}^{\mu_{[r,r']}(v)}\chi(\rho_m,\eta_m)\cdot1_{\{\rho_m\leq T\}}$
\end{enumerate}

We now fix the player II impulse control as  $\hat{v}=[\hat{\rho}_m,\hat{\eta}_m]_{m\geq1}\in \mathcal{V}$ and hence using (a) and (b) we find that $\nu(s)$ is now given by $\nu(s)=\sum_{m=1}^{\mu_{[t_0,s]}}\hat{\eta}_m\cdot1_{\{\rho_m\leq T\}}+\sum_{m=1}^{\mu_{[t_0,r]}}\xi_j\cdot1_{\{\tau_j\leq T\}}.  $\

As before, we employ the shorthand:
\begin{align}
    &Y^{y_0,\cdot}(s)\equiv (s,X^{t_0,x_0,\cdot}(t_0+s)), \quad y_0\equiv (t_0,x_0), \; \forall s\in [0,T-t_0], \\& \hat{Y}^{y_0,\cdot}(\tau)=Y^{y_0,\cdot}(\tau^{-} )+\Delta_N Y^{y_0,\cdot}(\tau), \quad \tau\in\mathcal{T},
\end{align}
where $\Delta_N Y(\tau) $ denotes a jump at time $\tau$ due to $\tilde{N}$.

Correspondingly, we adopt the following impulse response function $\hat{\Gamma}: \mathcal{T}\times S\times \mathcal{Z}\to  \mathcal{T}\times S$ acting on  $y'\equiv (\tau,x')\in \mathcal{T}\times S$ where $x'\equiv X^{t_0,x_0,\cdot}(t_0+\tau^-)$ is given by: 
\begin{align}
\hat{\Gamma}(y',\zeta)\equiv (\tau,\Gamma (x',\zeta))=(\tau,X^{t_0,x_0,\cdot} (\tau)),\quad \forall \xi\in\mathcal{Z},\; \forall\tau\in\mathcal{T} .
\end{align}

By It\={o}'s formula for c\`{a}dl\`{a}g semi-martingale (jump-diffusion) processes (see for example theorem II.33 of \cite{protter2005stochastic} in conjunction with Theorem 1.24 of \cite{oksendal2005applied}), we have that:
\begin{align}  
&\mathbb{E}[\phi (Y^{y_0  ,u, \hat{v} } (\tau_j  ))]- \mathbb{E}[\phi (\hat{Y}^{y_0  ,u, \hat{v} } (\tau_{j+1}^- ))]
\\&=- \mathbb{E}\left[\int_{\tau_j}^{\tau_{j+1}}\frac{\partial \phi }{\partial s}   +\mathcal{L}\phi (Y^{y_0  ,u, \hat{v} } (s ))ds  +\hspace{-5 mm} \sum_{\nu_{t_0,\tau_j}(\hat{v})< m\leq \nu_{t_0,\tau_{j+1}} (\hat{v})}\hspace{-6 mm}\Delta_\nu  \phi (Y^{y_0  ,u, \hat{v} } (\hat{\rho}_m  ))\right]. 	 \label{itoexpnthrm6.3}\end{align}

We note firstly that by definition of the intervention times $\{\tau_j\}_{j\in\mathbb{N}}$ we have that $\mu_{[\tau_j,\tau_{j+1}[}(u)=0$ since no player I interventions occur in the interval $[\tau_j,\tau_{j+1}[$. Hence, on the interval $[\tau_j,\tau_{j+1}[$ we have that $\Delta_\nu=\Delta_{\hat{\eta}}$ in particular, $\Delta_\nu  \phi =\Delta_{\hat{\eta}}  \phi$  so that $\Delta_\nu  \phi (Y^{y_0  ,u, \hat{v} } (\hat{\rho}_m  ))=\Delta_{\hat{\eta}}  \phi (Y^{y_0  ,u, \hat{v} } (\hat{\rho}_m  )) $.

Hence, by (\ref{verif_1_singular_cond}) we have that:
\begin{align}  
&\mathbb{E}[\phi (Y^{y_0  ,u, \hat{v} } (\tau_j  ))]- \mathbb{E}[\phi (\hat{Y}^{y_0  ,u, \hat{v} } (\tau_{j+1}^- ))]
\\&=- \mathbb{E}\left[\int_{\tau_j}^{\tau_{j+1}}\frac{\partial \phi }{\partial s}   +\mathcal{L}\phi (Y^{y_0  ,u, \hat{v} } (s ))ds  +\sum_{\nu_{t_0,\tau_j}(\hat{v})< m\leq \nu_{t_0,\tau_{j+1}} (\hat{v})}\hspace{-5 mm}\chi (\hat{\rho}_m  ,\hat{\eta}_m)\right].  \label{theorem6.3singtochicostconteqn} 
\end{align}
Summing both sides from $j=0$ to $j=k<\infty$ , we obtain the following:
\begin{align}
&\phi(y_0)+ \sum_{j=1}^{k} \mathbb{E}[\phi (\hat{Y}^{y_0  ,u, \hat{v} } (\tau_j  ))-\phi (\hat{Y}^{y_0  ,u, \hat{v} } (\tau_j ^-  ))]  -\mathbb{E}[\phi (\hat{Y}^{y_0  ,u, \hat{v} } (\tau_{k+1}^-  ))]
 	\\&\quad\leq \mathbb{E}\left[\int_{t_0}^{\tau_{k+1}}f(Y^{y_0,u,\hat{v}}(s))ds-\sum_{m\leq \nu_{t_0,\tau_{k+1}} (\hat{v})}\chi (\hat{\rho}_m  ,\hat{\eta}_m)    \right]. 	 
 	\label{theorem6.3aftersum}
 	\end{align}
Now by definition of the non-local intervention operator $\mathcal{M}_1$, we have that:
\begin{align*} \phi(Y^{y_0,u,\hat{v}} (\tau_j ))=\phi(\Gamma (\hat{Y}^{y_0,u,\hat{v}} (\tau_j-),\xi_j ))\geq \mathcal{M}_1 \phi(\hat{Y}^{y_0,u,\hat{v}} (\tau_j-))-c(\tau_j,\xi_j )\cdot 1_{\{\tau_j\leq \tau_S \}},\label{interventionineq1theorem6.3}	\end{align*}
(using the fact that $\inf_{z\in \mathcal{Z}}[\phi(\tau',\Gamma (X(\tau'^-),z))+c(\tau', z)\cdot 1_{\{\tau'\leq T \}}  ]=0$ whenever $\tau'>\tau_S$).

Hence,
\begin{align}
&\phi(Y^{y_0,u,\hat{v}} (\tau_j ))-\phi(\hat{Y}^{y_0,u,\hat{v}} (\tau_j-)) 
\\&\geq \mathcal{M}_1 \phi(\hat{Y}^{y_0,u,\hat{v}} (\tau_j-))-\phi(\hat{Y}^{y_0,u,\hat{v}} (\tau_j- ))-c(\tau_j,\xi_j )\cdot 1_{\{\tau_j\leq \tau_S \}},\label{interventionineq2theorem6.3} 
\end{align}
and by (\ref{vii_ch4_terminal_limit}) we readily observe that:
\begin{equation}
\phi(Y^{y_0,u,\hat{v}} (\tau_S ))-\phi(\hat{Y}^{y_0,u,\hat{v}} (\tau_S ))=0. 	
\end{equation}
After plugging (\ref{interventionineq2theorem6.3}) into (\ref{theorem6.3aftersum}) we obtain the following:
\begin{align*}
&\begin{aligned}
\phi (y_0)+ \sum_{j=1}^{k} \mathbb{E}[\mathcal{M}_1  \phi (\hat{Y}^{y_0  ,u, \hat{v} } (\tau_j^- ))-\phi (\hat{Y}^{y_0  ,u, \hat{v} } (\tau_j ^-  ))&-c (\tau_j  , \xi_j  ) \cdot 1_{\{\tau_j  \leq  \tau_S  \}}  ]
-\mathbb{E}[\phi (\hat{Y}^{y_0  ,u, \hat{v} } (\tau_{k+1}^-  ))]
\end{aligned}
 \\&\leq \mathbb{E}\left[\int_{t_0}^{\tau_{k+1}}f(Y^{y_0,u,\hat{v}} (s))ds-\hspace{-3 mm}\sum_{m\leq \nu_{t_0,\tau_{k+1}} (\hat{v})}\hspace{-4 mm}\chi (\hat{\rho}_m  ,\hat{\eta}_m)  \right]. 	
\label{prooftheorem6.3prelimit} 
\end{align*}
Hence,
\begin{align}
&\phi(y_0)+\sum_{j=1}^k\mathbb{E}[\mathcal{M}_1 \phi(\hat{Y}^{y_0,u,\hat{v}} (\tau_j^-))-\phi(\hat{Y}^{y_0,u,\hat{v}} (\tau_j^- ))]
\\&\begin{aligned}
\leq \mathbb{E}\Bigg[\int_{t_0}^{\tau_{k+1}}f(Y^{y_0,u,\hat{v}} (s))ds+\phi(\hat{Y}^{y_0,u,\hat{v}} (\tau_{k+1}^- ))&+\sum_{j\geq 1}c(\tau_j,\xi_j ) \cdot 1_{\{\tau_j\leq \tau_S \}} 
\\&-\sum_{m\leq \nu_{t_0,\tau_{k+1}}}\hspace{-2 mm}\chi (\hat{\rho}_m ,\hat{\eta}_m )    \Bigg], 	
\label{prooftheorem6.3prelimit3.1} 
\end{aligned}
\end{align}
Now $\lim_{k\to \infty }\sum_{j=1}^k\mathbb{E}[\mathcal{M}_1 \phi(\hat{Y}^{y_0,u,\hat{v}} (\tau_j-))-\phi(\hat{Y}^{y_0,u,\hat{v}} (\tau_j^- ))] =0$ since by (\ref{vi_ch4_terminal_limit}) we have that $\phi(\hat{Y}^{y_0,\cdot} (\tau_j ))-\phi(\hat{Y}^{y_0,\cdot} (\tau_j^- ))=0, \mathbb{P}-{\rm a.s}.$ when $\tau_j=\tau_S$; we can then deduce the statement by Lemma 3.10 in \cite{chen2013impulse} i.e. using the $\frac{1}{2}$-H\"{o}lder continuity of the non-local operator $\mathcal{M}_1$. Similarly we have by (\ref{vi_ch4_terminal_limit}) that $\phi(Y^{\cdot,} (s))\to G(Y^{\cdot,} (\tau_S ))$  as $s\to \tau_S^-$,  $\mathbb{P}-$a.s.
Hence, letting $k\to \infty$ in (\ref{prooftheorem6.3prelimit3.1}) gives:
\begin{align}\nonumber
\phi(y_0)\leq \mathbb{E}\Bigg[\int_{t_0}^{\tau_S}f(Y^{y_0,u,\hat{v}} (s))ds+\sum_{j\geq 1}c(\tau_j,\xi_j )  \cdot 1_{\{\tau_j\leq \tau_S \}}&-\sum_{n\geq 1}\chi (\hat{\rho}_n,\hat{\eta}_n )\cdot  1_{\{\hat{\rho}_n\leq \tau_S \}}
\\&+G(Y^{y_0,u,\hat{v}} (\tau_S))\cdot 1_{\{\tau_S<\infty\}}\Bigg]. 	 \label{prooftheorem6.3postlimit}\end{align}
Since this holds for all $u \in \mathcal{U}$, we have that:
\begin{align*}\nonumber
&\begin{aligned}\phi(y_0)\leq    \inf_{u \in \mathcal{U}} \mathbb{E}\Bigg[\int_{t_0}^{\tau_S}f (Y^{y_0,u, \hat{v} } (s ))ds  + \sum_{j\geq 1} c (\tau_j  , \xi_j  )   \cdot 1_{\{\tau_j  \leq  \tau_S  \}}  &- \sum_{n\geq 1} \chi  (\hat{\rho}_n  , \hat{\eta}_n  )\cdot  1_{\{\hat{\rho}_n  \leq  \tau_S  \}}  
\\&+G(Y^{y_0,u, \hat{v} } (\tau_S  ))\cdot 1_{\{\tau_S<\infty\}}\Bigg]  	 
\end{aligned}
\end{align*}
In particular, we have that:
\begin{align}&\nonumber
\begin{aligned}\phi(y_0)\leq    \sup_{v \in \mathcal{V}}\inf_{u \in \mathcal{U}}\mathbb{E}\Bigg[\int_{t_0}^{\tau_S}f(Y^{y_0,u,v} (s))ds+\sum_{j\geq 1}c(\tau_j,\xi_j )  \cdot 1_{\{\tau_j\leq \tau_S \}} -\sum_{n\geq 1}\chi (\rho_n,\eta_n )  \cdot 1_{\{\rho_n\leq \tau_S \}} &
\\
+G(Y^{y_0,x_0,u,v} (\tau_S))\cdot 1_{\{\tau_S<\infty\}}\Bigg]=V(y).&
\end{aligned}
\end{align}
Using an analogous argument, namely replacing $\hat{v}$ with $\hat{u}$ in (\ref{itoexpnthrm6.3}), then performing similar steps (using condition (\ref{iv_hjbi_p2_ineq})), we can similarly prove that:
\begin{align}\label{prooftheorem6.3finalP2}
\begin{aligned}\phi(y_0)
\geq \inf_{u \in \mathcal{U}}\sup_{v \in \mathcal{V}}\mathbb{E}\Bigg[\int_{t_0}^{\tau_S}f(Y^{y_0,u,v} (s))ds+\sum_{j\geq 1}c(\tau_j,\xi_j )\cdot  1_{\{\tau_j\leq \tau_S \}}-\sum_{m\geq 1}\chi (\rho_m,\eta_m ) \cdot 1_{\{\rho_m\leq \tau_S \}}& \\
+G(Y^{y_0,u,v} (\tau_S))\cdot 1_{\{\tau_S<\infty\}}\Bigg]&
\\=V(y).& 	 
\end{aligned}
\end{align}

Let us now fix the pair of controls $(\hat{u},\hat{v} )\in \mathcal{U}\times \mathcal{V}$, using the definition of $\Delta_z$ and by (\ref{verif_1_singular_cond}) we have that:
\begin{align}
0&=\Delta_z\phi(Y(\hat{\rho}_m)-\chi(\hat{\rho}_m,z)
\nonumber=\phi(\hat{\Gamma}(Y(\hat{\rho}_m^-,z)))-\phi(Y(\hat{\rho}_m^-))+\Delta_NY(\hat{\rho}_m)-\chi(\hat{\rho}_m,z)
\\&=\phi(\hat{\Gamma}(\hat{Y}(\hat{\rho}_m^-,z)))-\phi(Y(\hat{\rho}_m^-))-\chi(\hat{\rho}_m,z).
\label{marginalcostequation}
\end{align}

Now since (\ref{marginalcostequation}) holds for all $z\in\mathcal{Z}$, after applying the $\sup$ operator to both sides of (\ref{marginalcostequation}) we find that:
\begin{align}
0=\sup_{z\in\mathcal{Z}}[\phi(\hat{\Gamma}(\hat{Y}(\hat{\rho}_m^-),z))-\chi(\hat{\rho}_m,z)]-\phi(\hat{Y}(\hat{\rho}_m^-))
\nonumber=\mathcal{M}_2\phi(\hat{Y}(\hat{\rho}_m^-))-\phi(\hat{Y}(\hat{\rho}_m^-)),
\label{marginalcostequationtointervention}\end{align}
from which we immediately deduce the statement:
\begin{equation}
\mathcal{M}_2\phi(\hat{Y}(\hat{\rho}_m^-))=\phi(\hat{Y}(\hat{\rho}_m^-)).
\label{marginalcostequationtointervention2}\end{equation}
We now see that an immediate impulse intervention at $\hat{\rho}_m$ is indeed optimal for player II.
Using analogous arguments we can deduce that:
\begin{equation}
\mathcal{M}_1\phi(\hat{Y}(\hat{\tau}_j^-))=\phi(Y(\hat{\tau}_j^-)).
\label{marginalcostequationtointerventionP2}\end{equation}
We hereafter straightforwardly observe using (\ref{v_ch4_continuity_equation_equality}) and (T4) we find the following \textit{equality}:
\begin{align}\nonumber
\phi(y_0)=\mathbb{E}\Bigg[\int_{t_0}^{\tau_S}f(Y^{y_0,\hat{u},\hat{v}} (s))ds+\sum_{j\geq 1}c(\hat{\tau}_j,\hat{\xi}_j ) \cdot 1_{\{\hat{\tau}_j\leq \tau_S \}} &-\sum_{m\geq 1}\chi (\hat{\rho}_m,\hat{\eta}_m )  \cdot 1_{\{\hat{\rho}_m\leq \tau_S \}}
\\&+G(Y^{y_0,\hat{u},\hat{v}} (\tau_S ))\cdot 1_{\{\tau_S<\infty\}}\Bigg].	 \label{theorem6.3valuefunctionequalityfinal}
\end{align}
Hence, we can deduce the following statement:
\begin{equation} \sup_{v \in \mathcal{V}}\inf_{u \in \mathcal{U}}J^{u,v}[y]\geq \phi(y)=J^{\hat{u},\hat{v}}[y]\geq \inf_{u \in \mathcal{U}}\sup_{v \in \mathcal{V}}J^{u,v}[y],\quad \forall  y\in [0,T]\times S. \label{pinchingvaluefunctionequation6.3}	
\end{equation}
Now, since $\inf_{u \in \mathcal{U}}\sup_{v \in \mathcal{V}}J^{u,v}[y]\geq \sup_{v \in \mathcal{V}}\inf_{u \in \mathcal{U}}J^{u,v}[y]$($=V^{+}(y)$) it then follows that:
\begin{equation} 
V(y)=\phi(y)=J^{\hat{u},\hat{v}}[y],\qquad \forall  y\in [0,T]\times S,\label{phieqV} 
\end{equation}
after which we deduce the thesis.
\end{refproof}

\section{A HJBI Equation for Non-Zero-Sum Stochastic Differential Games with Impulse Controls} \label{section_verification_non_zero_sum}

In this section, we now extend the results to non zero-sum stochastic differential games. Expectedly, proving existence results of Nash equilibria for stochastic differential games that involve impulse controls relies on a similar set of arguments as those constructed in the continuous control case. Indeed, Nash equilibria existence and characterisation results have been established for differential games in which continuous controls were used see e.g. \cite{KononenkoNonantagonisticDifferentialGames,KleimenovNonantagonisticDifferentialGames} and in the stochastic case in \cite{buckdahn2004nash} in which a method of generalising Folk Theorems\footnote{Folk Theorems are a set of fundamental results within repeated games that state that any feasible payoff for which the player is weakly better off than their minmax payoff (i.e. individual rational payoff) can in fact be supported in subgame perfect Nash equilibrium when the players are
sufficiently patient that is, whenever the players of the game are sufficiently patient then the repeated game can allow any outcome in the average payoff sense.} in classical game theory (i.e. deterministic repeated games) was used.

We firstly prove a non zero-sum verification theorem for the game in which both players use impulse controls to modify the state process.

In order to describe a non-zero sum game, we now consider a game in which each player has their own individual payoff function. The  payoff functions for player I and player II, $J_1$ and $J_2$ respectively, are given by the following:
\begin{align}
&\begin{aligned}
J_1^{(\tilde{u},\tilde{v} )} [t_0,x_0]=\mathbb{E}\Bigg[\int_{t_0}^{\tau_S}f_1(s,X^{t_0,x_0,\tilde{u},\tilde{v}} (s))ds&-\sum_{j\geq 1}c_1 (\tilde{\tau}_j,\tilde{\xi }_j ) \cdot 1_{\{\tau_j\leq \tau_S \}}
\\&+G_1 (\tau_S,X^{t_0,x_0,\tilde{u},\tilde{v}} (\tau_S ))\cdot 1_{\{\tau_S<\infty\}}\Bigg], \label{NEJ1}	
\end{aligned}\\&\begin{aligned}
J_2^{(\hat{u},\tilde{v} )}  [t_0,x_0]=\mathbb{E}\Bigg[\int_{t_0}^{\tau_S}f_2 (s,X^{t_0,x_0,\tilde{u},\tilde{v}} (s))ds&-\sum_{m\geq 1}c_2 (\tilde{\rho}_m,\tilde{\eta}_m ) \cdot 1_{\{\rho_m\leq \tau_S \}}
\\&+G_2  (\tau_S,X^{t_0,x_0,\tilde{u},\tilde{v}} (\tau_S ))\cdot 1_{\{\tau_S<\infty\}} \Bigg],
\label{NEJ2}	
\end{aligned}\\&\nonumber
\qquad\qquad \qquad \qquad \qquad\qquad\qquad\qquad\qquad\qquad\qquad\qquad\qquad  \forall (t_0,x_0)\in  [0,T]\times S,
\end{align} 
where $\tilde{u}=[\tilde{\tau}_j,\tilde{\xi}_j]_{j\geq 1}$ and $\tilde{v}=[\tilde{\rho}_m,\tilde{\eta}_m]_{m\geq 1}$ are admissible controls for player I and player II respectively. 

We note that the function $J_1^{(\tilde{u},\hat{v})} (t,x) $ (resp., $J_2^{(\hat{u},\tilde{v} )} (t,x)) $ defines the payoff received by the player I (resp., player II) when it uses the control $\tilde{u}\in \mathcal{U}$ (resp., $\tilde{v}\in \mathcal{V}$) and player II (resp., player I) uses the control $\hat{v}\in \mathcal{V}$ (resp. $\hat{u}\in \mathcal{U}$) given some initial point $(t,x)\in [0,T]\times S$.

The following are key objects in the analysis of impulse control models:

\begin{definition}\label{Intervention_operator_two_players} Let $\tau\in\mathcal{T}$. For $i\in\{1,2\}$, we define the [non-local] Player $i$-intervention operator $\mathcal{M}_i:\mathcal{H}\to \mathcal{H}$ acting at a state $X(\tau)$ by the following expression:
\begin{equation} 
\mathcal{M}_i \phi(\tau,X(\tau)):=\sup_{z\in \mathcal{Z}}[\phi(\tau,\Gamma (X(\tau^-),z))-c_i (\tau, z)\cdot 1_{\{\tau\leq T\}}  ], \label{NEinterventionop}\end{equation}
where $\Gamma : S\times \mathbb{R}\to  S$ is the impulse response function.
\end{definition}
\begin{definition}[Nash Equilibrium for Non-Zero-Sum Games with Impulse Control]\label{Nash_Equilibrium_for_Non-Zero-Sum_Games_with_Impulse_Control}

We say that a pair $(\hat{u},\hat{v} )\in \mathcal{U}\times \mathcal{V}$ is a Nash equilibrium of the stochastic differential game with impulse controls $\hat{u}=[\hat{\tau}_j,\hat{\xi}_j ]_{j\in \mathbb{N}}\in \mathcal{U}, \hat{v}=[\hat{\rho}_m,\hat{\eta}_m ]_{m\in \mathbb{N}}\in \mathcal{V}$ if the following statements hold:
 \renewcommand{\theenumi}{\roman{enumi}}
 \begin{enumerate}
\item 
	$J_1^{(u,\hat{v} )} [t,x]\geq J_1^{(\hat{u},\hat{v}) } [t,x] \quad \forall u \in \mathcal{U},\forall  (t,x)\in [0,T]\times S, \label{p_1_NE_cond_ch4}$
\item
	$J_2^{(\hat{u},v)} [t,x]\geq J_2^{(\hat{u},\hat{v} )} [t,x] \quad \forall v \in \mathcal{V},\forall  (t,x)\in [0,T]\times S. \label{p_2_NE_cond_ch4}$
\end{enumerate}
\end{definition}
Condition \eqref{p_1_NE_cond_ch4} states that given some fixed player II control policy $\hat{v}\in \mathcal{V}$, player I cannot profitably deviate from playing the control policy $\hat{u}$. Analogously, condition \eqref{p_2_NE_cond_ch4} is the equivalent statement given the player I's control policy is fixed as $\hat{u}$, player II cannot profitably deviate from $\hat{v}$. We therefore see that $ (\hat{u},\hat{v})$ is an equilibrium in the sense of a Nash equilibrium since neither player has an incentive to deviate whenever their opponent plays their equilibrium policy.

We now generalise our the zero-sum result of Theorem \ref{Verification_theorem_for_Zero-Sum_Games with Impulse Control} to cover non zero-sum payoff structures with the use of a Nash Equilibrium solution concept. 
\begin{theorem}[Verification theorem for Non-Zero-Sum Games with Impulse Control]\label{Verification_theorem_for_Non-Zero-Sum_Games_with_Impulse_Control}

Let us suppose that the value of the game exists and that there exists functions $\phi_i\in \mathcal{C}^{1,2} ([0,T],S)\cap\mathcal{C}([0,T],\bar{S}),i\in \{1,2\}$ s.th. $\phi_i$ satisfy technical conditions (T1) - (T4) and the following conditions $\forall s\in[0,T]$:
 \begin{enumerate}[label={(\roman*')}]	
\item	$\phi_i\geq \mathcal{M}_i \phi_i$ on $S$ and the regions $D_i$ are defined by:\\ $\noindent D_i=\{x\in S;\phi_i (\cdot,x)>\mathcal{M}_i \phi_i (\cdot,x)\},   i\in \{1,2\}$
where we refer to $D_1$ (resp., $D_2$) as the player I (resp., player II) continuation region.
	\item	$\frac{\partial \phi_1}{\partial s}+\mathcal{L}\phi_1 (\cdot,X^{\cdot,u ,\hat{v}} (\cdot))+f_1 (\cdot X^{\cdot,u ,\hat{v}} (\cdot))\\\geq \frac{\partial \phi_1}{\partial s}+\mathcal{L}\phi_1 (\cdot,X^{\cdot,\hat{u} ,\hat{v}} (\cdot))+f_1 (\cdot,X^{\cdot,\hat{u} ,\hat{v}} (\cdot))\geq 0, \hspace{3 mm}  \hspace{1 mm}\forall  u \in \mathcal{U}$ on $S\backslash{D_1}$.
	\item $\frac{\partial \phi_2}{\partial s}+\mathcal{L}\phi_2 (\cdot,X^{\cdot,\hat{u} ,v  } (\cdot))+f_2 (\cdot,X^{\cdot,\hat{u} ,v  } (\cdot))\\\geq \frac{\partial \phi_2}{\partial s}+\mathcal{L}\phi_2 (\cdot,X^{\cdot,\hat{u} ,\hat{v}} (\cdot))+f_2 (\cdot,X^{\cdot,\hat{u} ,\hat{v}} (\cdot))\geq 0, \hspace{3 mm}  \hspace{1 mm}\forall  v \in \mathcal{V}$ on $S\backslash{D_2}$.
	\item \label{case_2_equality_continuity}	$\frac{\partial \phi_i}{\partial s}+\mathcal{L}\phi_i (\cdot,X^{\cdot,\hat{u},\hat{v} } (\cdot))+f_i (\cdot,X^{\cdot,\hat{u}_{[t_0,s]} ,\hat{v}_{[t_0,s]}} (\cdot))=0$ on $D_1\cap D_2$.
	
	\item	$\hat{\xi}_k\in \arg\hspace{-0.55 mm}\sup_{z\in \mathcal{Z}}\{\phi_i (\tau_k,\Gamma (x,z))-c(\tau_k,z)\}$ is a Borel Measurable selection $\forall x \in S,\;\tau_k\in\mathcal{T}$. Similarly, $\hat{\eta}_j\in \arg\hspace{-0.55 mm}\sup_{z\in \mathcal{Z}}\{\phi_i (\rho_j,\Gamma (x,z))-\chi (\rho_j,z)\}$ is a Borel Measurable selection $\forall x \in S,\;\forall \rho_j\in\mathcal{T}$.
\end{enumerate}

Put $\hat{\tau}_0\equiv t_0$ and define $\hat{u}:=[\hat{\tau}_j,\hat{\xi}_j ]_{j\in \mathbb{N}}$ inductively by $\hat{\tau}_{j+1}=\inf\{s>\tau_j;\;X^{\hat{u}_{[t_0,s]} ,v} (\cdot)\notin D_1 \}\wedge \tau_S,   \hspace{1 mm}\forall  v \in \mathcal{V}$ . Similarly, put $\hat{\rho}_0\equiv t_0$ and define $\hat{v}:=[\hat{\rho}_m,\hat{\eta}_m ]_{m\in \mathbb{N}}$ inductively by $\hat{\rho}_{m+1}=\inf\{s>\rho_m;X^{u,\hat{v}_{[t_0,s]}} (\cdot)\notin D_2 \}\wedge \tau_S,  \hspace{1 mm}\forall  u \in \mathcal{U}$.  

Then $(\hat{u},\hat{v})$ is a Nash equilibrium for the game, that is to say the following statements hold:
\begin{equation}
\phi_1 (t,x)=\sup_{u \in \mathcal{U}}J_1^{(u,\hat{v} )} [t,x]=J_1^{(\hat{u},\hat{v} )} [t,x],	\end{equation}
and
\begin{equation} \phi_2 (t,x)=\sup_{v \in \mathcal{V}}J_2^{(\hat{u},v)} [t,x]=J_2^{\hat{u},\hat{v} } [t,x].	 \end{equation}
\hfill$\forall  (t,x)\in [0,T]\times S$.
\end{theorem}

% As in Theorem \ref{Verification_theorem_for_Zero-Sum_Games with Impulse Control}, we note that conditions \ref{verif_theorem_c-s_NON_zero_item_ii'_intervention_stopping_inequalities} - \ref{verif_theorem_c-s_NON_zero_item_iv'_HJBI_equation} of Theorem \ref{Verification_theorem_for_Non-Zero-Sum_Games_with_Impulse_Control} follow directly from QVI conditions which can be  motivated by a heuristic analysis similar to that of the zero-sum case. 
The condition $\phi_i\in \mathcal{C}^{1,2} ([0,T],S)\cap\mathcal{C}([0,T],\bar{S}),\; i\in \{1,2\}$ is used to allow for the integro-differential operator $\mathcal{L}$ in (ii') - (iv') to be applied in addition to permitting an application of Dynkin's formula which is central to the proof of the theorem.

The proof of Theorem \ref{Verification_theorem_for_Non-Zero-Sum_Games_with_Impulse_Control} follows a similar path to that of Theorem \ref{Verification_theorem_for_Zero-Sum_Games with Impulse Control}, we therefore defer the proof of the theorem to the appendix. 

% In the following, we apply Theorem \ref{Verification_theorem_for_Non-Zero-Sum_Games_with_Impulse_Control} to solve the duopoly investment problem and in doing so, provide an example of a function that satisfies the properties of the theorem. 

The proof of Theorem \ref{Verification_theorem_for_Non-Zero-Sum_Games_with_Impulse_Control} follows a similar path to that of Theorem \ref{Verification_theorem_for_Zero-Sum_Games with Impulse Control}. We therefore defer the proof of the theorem to the appendix. In the next section we apply Theorem \ref{Verification_theorem_for_Non-Zero-Sum_Games_with_Impulse_Control} to solve the duopoly investment problem in Section \ref{section_duopoly} and in doing so, provide an example of a function that satisifies the properties of the theorem. \

Before doing so, in analogy to Corollary \ref{Corollary 6.2.}, we give the following result which follows directly from Theorem \ref{Verification_theorem_for_Non-Zero-Sum_Games_with_Impulse_Control}:

\begin{corollary}\label{Corollary 7.4.} When each player plays their equilibrium control, the sample space splits into three regions that represent a region in which player I intervenes in $I_1$, a region in which player II intervenes $I_2$, and a region in which no action is taken by either player I or player II; moreover the three regions are characterised by the following expressions for $j\in \{1,2\}$:
\begin{align*}
I_j&=\{(t,x)\in [0,T]\times S: V_j (t,x)=\mathcal{M}_j V_j (t,x),\;\mathcal{L}V_j (t,x)+f_j (t,x)\geq 0\}  .\\
I_3&=\{(t,x)\in [0,T]\times S: V_j (t,x)\geq \mathcal{M}_j V_j (t,x);\; \mathcal{L}V_j (t,x)+f_j (x)=0\}.
\end{align*}
\end{corollary}
\section{Examples} \label{section_examples}
In order to demonstrate the workings of the theorems, we give some example calculations.  \begin{example}

Consider a system with passive dynamics that are described by the following stochastic process:
\begin{align}
dX(r)=\alpha dr+\beta dB(r),\qquad \forall r\in [0,T], \label{example_calc_ch4_zs_state_process}
\end{align}
where $\alpha,\beta\in\mathbb{R}_{>0}$ are fixed constants, $B(r)$ is a $1$-dimensional Brownian motion and $T\in\mathbb{R}_{>0}$ is some finite time horizon. The process $X$ in (\ref{example_calc_ch4_zs_state_process}) is known as \textit{Brownian motion with drift} and models and number of processes in finance such as insurance claim processes and  risk-neutral price-processes in options pricing, for example \cite{david2005minimizing, pechtl1999some}. 

The state process $X$ is modified by two controllers, player I, that exercises an impulse control policy $u=[\tau_j,\xi_j]\in\mathcal{U}$ and player II that exercises an impulse control policy $v=[\rho_m,\eta_m]\in\mathcal{V}$. The controlled state process evolves according to the following expression:
\begin{align}\nonumber
X(t)=x_0+ \alpha\int_0^{t\wedge\tau_S} ds+\beta\int_0^{t\wedge\tau_S} dB(s)&-\sum_{j\geq 1}(\kappa_1+(1+\lambda_1)\xi_j)\cdot 1_{\{\tau_j\leq t\wedge\tau_S\}}
\\&
\begin{aligned}-\sum_{m\geq 1}(\kappa_2+(1+\lambda_2)\eta_m)\cdot 1_{\{\rho_m\leq t\wedge\tau_S\}},&\label{example_calc_ch4_zs_control_state_process}
\\X(0)\equiv x_0,\nonumber
\\\forall t\in [0,T],\mathbb{P}-{\rm a.s.},&\nonumber
\end{aligned}
\end{align}
where $\tau_S:=\inf\{s>0:X(s)\leq 0\}$ and the constants $\kappa_i>0$ and $\lambda_i>0$ are the fixed part and the proportional part of the transaction cost incurred by player $i\in\{1,2\}$ for each intervention (resp.).

Player I seeks to choose an admissible impulse control $u=[\tau_j,\xi_j]$ that maximises its reward $J$ where $\{\tau_j\}_{\{j\geq 1\}}$ are player I intervention times and each $\xi_{j\geq 1}\in\mathcal{Z}$ is a player I impulse intervention. Player II seeks to choose an admissible impulse control $v=[\rho_m,\xi_m]$ that minimises the same quantity $J$ where $\{\rho_m\}_{\{m\geq 1\}}$ are player II  intervention times and each $\eta_{m\geq 1}\in\mathcal{Z}$ is a player II impulse intervention. 

The function $J$ is given by the following expression:
\begin{align}
J^{u,v}[t,x]=\mathbb{E}\left[\sum_{j\geq 1}e^{-\delta\tau_j}\xi_j\cdot 1_{\{t<\tau_j\leq \tau_S\}}-\sum_{m\geq 1}e^{-\delta\rho_m}\eta_m\cdot 1_{\{t<\rho_m\leq \tau_S\}} \right] , \qquad \forall (t,x)\in[0,T]\times\mathbb{R},
\end{align}
where $\delta \in ]0,1]$ is common discount factor.

An example of a setting for this game is an interaction between two players  that consume a common rivalrous and exhaustible good (e.g. public funds, extractable resources, labour supply etc.) which is vulnerable to stochastic extinction. For each act of consumption, each player incurs both a fixed cost and a proportional cost.

The problem is to find a function $\phi\in\mathcal{C}^{1,2}([0,T],\mathbb{R})$ s.th.
\begin{align}
\sup_{u}\inf_{v}J^{u,v}(s,x)=\inf_{v}\sup_{u}J^{u,v}(s,x)=\phi(s,x),\quad \forall (s,x)\in[0,T]\times\mathbb{R}.
\end{align}
% which by Theorem  \ref{Verification_theorem_for_Zero-Sum_Games with Impulse Control} is a candidate for the value function of the game.

We recognise this as a zero-sum stochastic game in which both players use impulse controls. We therefore seek to apply Theorem \ref{Verification_theorem_for_Zero-Sum_Games with Impulse Control} to compute the equilibrium controls. 

Firstly, we observe that by (\ref{example_calc_ch4_zs_state_process}) and using (\ref{generator}), the generator $\mathcal{L}$ for the process $X$ is given by: 
\begin{align}
\mathcal{L}\psi(s,x)=\frac{\partial \psi}{\partial s}(s,x)+\alpha  \frac{\partial \psi}{\partial x}(s,x)+\frac{1}{2}\beta^2\frac{\partial^2 \psi}{\partial x^2}(s,x),    \label{example_ch4_generator}
\end{align}
for some test function $\psi\in\mathcal{C}^{1,2}([0,T],\mathbb{R})$.

We now wish to derive the functional form of the function $\phi$. Applying (\ref{v_ch4_continuity_equation_equality}) of Theorem \ref{Verification_theorem_for_Zero-Sum_Games with Impulse Control} leads to the \textit{heat equation} $\mathcal{L}\phi=0$ (here, $f\equiv 0$ in Theorem \ref{Verification_theorem_for_Zero-Sum_Games with Impulse Control}). Following this, we make the following ansatz for the function $\phi$, $\phi(s,x)=e^{-\delta s}\psi(x),\psi(x):=ae^{bx}$ for some as yet, undetermined constants $a,b\in\mathbb{R}$.

Plugging the ansatz for the function $\phi$ and using (\ref{v_ch4_continuity_equation_equality}) of Theorem \ref{Verification_theorem_for_Zero-Sum_Games with Impulse Control} into (\ref{example_ch4_generator}) immediately gives:
\begin{align}
-\delta+\alpha b +\frac{1}{2}\beta^2b^2 =0.    
\end{align}
After some manipulation, we deduce that there exist two solutions for $b$ which we denote by $b_1$ and $b_2$ s.th. $b_1>b_2$ with $b_1>0$ and $|b_2|>0$ which are given by the following:
\begin{align}
b_1=\frac{1}{\beta^2}\sqrt{\alpha^2+2\beta^2\delta}-\frac{\alpha}{\beta^2}, \qquad b_2=-\frac{\alpha}{\beta^2}-\frac{1}{\beta^2}\sqrt{\alpha^2+2\beta^2\delta}.
\label{b_values_ch4_example}
\end{align}
We now apply the HJBI equation (\ref{v_ch4_continuity_equation_equality}) of Theorem \ref{Verification_theorem_for_Zero-Sum_Games with Impulse Control} to characterise the function $\phi$ on the region $D_1\cap D_2$. Following our ansatz, we  observe, using (\ref{v_ch4_continuity_equation_equality}), the following expression for the function $\phi$: 
\begin{align}
\phi(s,x)&=e^{-\delta s}\psi(x),\qquad &&\forall (s,x)\in[0,T]\times D_1\cap D_2,\\
\psi(x)&=(a_1e^{b_1x}+a_2e^{b_2x}), \qquad &&\forall x\in D_1\cap D_2, \label{ch_4_phi_a}
\end{align}
where $a_1$ and $a_2$ are constants that are yet to be determined and $D_1$ and $D_2$ are the continuation regions for player I and player II respectively.

In order to determine the constants $a_1$ and $a_2$, we firstly observe that $\phi(\cdot,0)=0$. This implies that $a_1=-a_2:=a>0$. We now deduce that the function $\psi$ is given by the following expression:
\begin{align}
\psi(x)=a(e^{b_1x}-e^{b_2x}),\qquad \forall x\in D_1\cap D_2.      
\end{align}
In order to characterise the function over the entire state space and find the value $a$, using conditions (\ref{ii_ch4_intervention_ineq}) - (\ref{verif_1_singular_cond}) of Theorem \ref{Verification_theorem_for_Zero-Sum_Games with Impulse Control}, we study the behaviour of the function $\phi$ given each of the players' equilibrium controls.

Firstly, we consider the player I impulse control problem. In particular, we seek to obtain conditions on the impulse intervention applied when $\mathcal{M}_1\phi=\phi$. To this end, let us firstly conjecture that the player I continuation region $D_1$ takes the following form:
\begin{equation}
    D_1=\{x\in\mathbb{R}; 0<x<\tilde{x}\},\label{ch4_example_cont_region_d1}
\end{equation}
for some constant $\tilde{x}$ which we shall later determine.

Our first task is to determine the optimal value of the impulse intervention. We now define the following two functions which will be of immediate relevance:
\begin{align}
&\psi_0(x):=a(e^{b_1x}-e^{b_2x}),\\
&h(\xi):=\psi(x-\kappa_1-(1+\lambda_1)\xi)+\xi,\\&\qquad\qquad\qquad\qquad\qquad\qquad\qquad\qquad\qquad\qquad
\forall x\in \mathbb{R},\forall \xi\in\mathcal{Z}.\nonumber
\end{align}
In order to determine the value $\hat{\xi}$ that maximises $\Gamma(x({\tau}^- ),\xi)$ at the point of intervention, we investigate  the first order condition on $h$ i.e. $0=h'(\xi)$. This implies the following:
\begin{align}
\psi'(\tilde{x}-\kappa_1-(1+\lambda_1)\xi)&=\frac{1}{1+\lambda_1}.\label{ch_4_example_impulse_1st_condition}
\end{align}
Using the expression for $\psi$ (\ref{ch_4_phi_a}), we also observe the following:
\begin{align}
\psi'_0(x)= b_1e^{b_1x}&-b_2e^{b_2x}>0,\qquad&&\forall x\in\mathbb{R},  \label{ch_4_example_FOC}
\\
\psi''_0(x)=b_1^2e^{b_1x}&-b_2^2e^{b_2x}<0,&&\forall x<x^{\#}:=\frac{2}{b_1-b_2}\ln{\Bigg[\frac{|b_2|}{b_1}\Bigg]},
\label{ch_4_example_SOC}
\end{align}
from which we deduce the existence of two points $x^{\star},x_{\star}$ for which the condition $\psi'_0(x)=(1+\lambda_1)^{-1}$ holds. W.l.o.g. we assume $x^{\star}>x_{\star}$. Now by (\ref{ii_ch4_intervention_ineq}) of Theorem \ref{Verification_theorem_for_Zero-Sum_Games with Impulse Control} we require that $\phi(\cdot,x)=\mathcal{M}_1\phi(\cdot,x)$ whenever $x\geq \tilde{x}$ (c.f. $D_1$ in equation (\ref{ch4_example_cont_region_d1})) for which we have $e^{-\delta t}\psi_0(x)=\mathcal{M}_1\phi(t,x)$ whenever $x\geq \tilde{x}$, hence we find that:
\begin{equation}
    \psi(x)=\psi_0(x_{\star})+\hat{\xi}(x), \qquad \forall x\geq \tilde{x},
\end{equation}
where $x-\kappa_1-(1+\lambda_1)\hat{\xi}(x)=x_{\star}$ from which we readily find that the optimal player I impulse intervention value is given by:
\begin{equation}
    \hat{\xi}(x)=\frac{x-x_{\star}-\kappa_1}{1+\lambda_1},\qquad \forall x\geq \tilde{x}.\label{ch_4_example_optim_intervention}
\end{equation}
Having determined the optimal impulse intervention and constructed the form of the continuation region for Player I, we can derive the optimal impulse intervention for player II straightforwardly by analogous arguments from which we find that the continuation region for player II takes the form:
\begin{equation}
    D_2=\{x\in\mathbb{R}; 0<x<\bar{x}\},\label{ch4_example_cont_region_d2}
\end{equation}
and the optimal player II impulse intervention value is given by
\begin{equation}
    \hat{\eta}(x)=\frac{x_{\#}-x-\kappa_2}{1+\lambda_2},\qquad \forall x\geq \bar{x}.\label{ch_4_example_optim_intervention_pII}
\end{equation}

Putting the above facts together yields the following characterisation of the function $\psi$:
\begin{align}
\psi(x) =
\begin{cases}
\begin{aligned}
&a(e^{b_1x}-e^{b_2x}),&& &\forall x\in D_1\cap D_2,
\\&a(e^{b_1x_{\#}}-e^{b_2x_{\#}})&&\hspace{-2.8 mm}+\frac{x_{\#}-x-\kappa_2}{1+\lambda_2}, &\forall x\notin  D_2,
\\&a(e^{b_1x_{\star}}-e^{b_2x_{\star}})&&\hspace{-2.8 mm}+\frac{x-x_{\star}-\kappa_1}{1+\lambda_1}, &\forall x\notin D_1,
\end{aligned}
\end{cases}
\end{align}
where the constants $b_1$ and $b_2$ are specified in equation (\ref{b_values_ch4_example}).

Using the facts above and, invoking the high contact principle\footnote{Recall that the high contact principle is a condition that asserts the continuity of the value function at the boundary of the continuation region.}, we are now in a position to determine the value of the constants $a,\hat{x}$ and $\tilde{x}$. 

We now apply the high contact principle to find the boundary of the continuation region $D_2$. Indeed, continuity at $\bar{x}$ leads to the following:
\begin{align}
\psi(\bar{x})=\psi_0(x_{\#})+\hat{\eta}(\bar{x})  \implies a(e^{b_1\bar{x}}-e^{b_2\bar{x}})=a(e^{b_1x_{\#}}-e^{b_2x_{\#}})+\frac{x_{\#}-\bar{x}-\kappa_2}{1+\lambda_2}, 
\end{align}
from which we find that $\bar{x}$ is the solution to the following equation:
\begin{align}
    m_2(\bar{x})&=0,\label{ch_4_example_p2_m_function}
    \end{align} 
where the function $m_2$ is given by:
\begin{align}
        m_2(x)&=x-a(1+\lambda_2)(e^{b_1x}-e^{b_2x}-e^{b_1x_{\#}}+e^{b_2x_{\#}} )+x_{\#}-\kappa_2. \label{ch_4_example_p2_cont_region_value}
\end{align}

Lastly, we reapply the high contact principle to find the boundary of the continuation region $D_1$. Indeed, continuity at $\tilde{x}$ leads to the following relationship:
\begin{align}
\psi(\tilde{x})=\psi_0(x_{\star})+\hat{\xi}(\tilde{x})  \implies a(e^{b_1\tilde{x}}-e^{b_2\tilde{x}})=a(e^{b_1x_{\star}}-e^{b_2x_{\star}})+\frac{\tilde{x}-x_{\star}-\kappa_1}{1+\lambda_1},
\end{align}
from which we find that $\tilde{x}$ is the solution to the following equation:
\begin{align}
    m_1(\tilde{x})=0,\label{ch_4_example_p1_m_function}
\end{align}
where the function $m$ is given by:
\begin{align}
        m_1(x)=x-a(1+\lambda_1)(e^{b_1x}-e^{b_2x}-e^{b_1{x}_{\star}}+e^{b_2{x}_{\star}} )-{x}_{\star}-\kappa_1. \label{ch_4_example_p1_cont_region_value}
\end{align}

Equations (\ref{ch_4_example_p2_cont_region_value}) and (\ref{ch_4_example_p1_m_function}) are difficult to solve analytically for the general case but can however, be straightforwardly solved numerically using a root-finding algorithm.

To summarise, the solution is as follows: whenever $X\in D_1\cap D_2$ neither player intervenes. Player I performs an impulse intervention of size $\hat{\xi}$ given by (\ref{ch_4_example_optim_intervention}) whenever the process reaches the value $\tilde{x}$ and player II performs an impulse intervention of size $\hat{\eta}$ given by (\ref{ch_4_example_optim_intervention_pII}) whenever the process reaches the value $\bar{x}$. The value function for the problem is $\phi(s,x)= e^{-\delta s}\psi(x),\; \forall (s,x)\in[0,T]\in\mathbb{R}$, where is $\psi$ given by:
\begin{align}
\psi(x) =
\begin{cases}
\begin{aligned}
&a(e^{b_1x}-e^{b_2x}),&& &\forall x\in D_1\cap D_2,
\\&a(e^{b_1x_{\#}}-e^{b_2x_{\#}})&&\hspace{-2.8 mm}+\frac{x_{\#}-x-\kappa_2}{1+\lambda_2}, &\forall x\notin  D_2,
\\&a(e^{b_1x_{\star}}-e^{b_2x_{\star}})&&\hspace{-2.8 mm}+\frac{x-x_{\star}-\kappa_1}{1+\lambda_1}, &\forall x\notin D_1,
\end{aligned}
\end{cases}
\end{align}
and where the player I and player II continuation regions are given by:
\begin{align}
    &D_1=\{x\in\mathbb{R}\;; 0<x<\tilde{x}\},
    \\ &D_2=\{x\in\mathbb{R}\;; 0<x<\bar{x}\},
\end{align}
where the constants $\bar{x}$ and $\tilde{x}$ are determined by (\ref{ch_4_example_p2_cont_region_value}) and (\ref{ch_4_example_p1_cont_region_value}) respectively and the constants $b_1,b_2$ are given by (\ref{b_values_ch4_example}).
\end{example}
\begin{example}
\subsubsection*{The Duopoly Investment Problem Revisited}
We apply our results to prove Theorem \ref{Theorem 4.3.}. Let us denote by $Y$ the process $Y(s)=(s+t_0,X_1(s),X_2(s))$, where $X_1 :\mathbb{R}_{>0}\times \Omega \to \mathbb{R},\; X_2 :\mathbb{R}_{>0}\times \Omega \to \mathbb{R}$ are processes which represent the market share processes for Firm 1 and Firm 2  respectively and whose evolution is described by (\ref{case3stateprocessFirm1}) - (\ref{case3stateprocessFirm2}). We wish to fully characterise the optimal investment strategies for each firm, in order to do this we apply Theorem \ref{Verification_theorem_for_Non-Zero-Sum_Games_with_Impulse_Control}. We restrict ourselves to the case when $\theta_{ij}(s)\equiv\bar{\theta}_{ij}\in\mathbb{R}\backslash\{0\}$ and $\sigma_{ij}(s)\equiv\bar{\sigma}_{ij}\in\mathbb{R}\backslash\{0\}$.

\begin{refproof}[Proof of Theorem \ref{Theorem 4.3.}]
Given an admissible Firm 1 (resp., Firm 2) investment policy $u=[\tau_j,\xi_j]_{j\in\mathbb{N}}\in\mathcal{U}$ (resp., $v=[\rho_m,\eta_m]_{m\in\mathbb{N}}\in\mathcal{V}$) we note that the following identities hold:
\begin{align} 
X_1(\tau_j)&=\Gamma(X_1(\tau_j-)+\Delta_NX_1(\tau_j),\xi_j)=\hat{X}_1(\tau_j)+\xi_j, \label{impulserespeqnP1}
\\ X_2(\rho_m)&=\Gamma(X_2(\rho_m-)+\Delta_NX_2(\rho_m),\eta_m)=\hat{X}_2(\rho_m)+\eta_m.
\label{impulserespeqnP2}\end{align}

The Firm 1 and Firm 2 investment intervention operator (acting on some function $\psi:[0,T]\times S\to\mathbb{R}$) are given by the following expressions:
\begin{align}
\mathcal{M}_1\psi(\tau^-,x)&=\sup_{\xi\in\mathcal{Z}}\{\psi(\tau,x+\xi)-(\lambda_1\xi+\kappa_1)\},\label{interventionP1duop}
\\\mathcal{M}_2\psi(\tau^-,x)&=\sup_{\eta\in\mathcal{Z}}\{\psi(\tau,x+\eta)-(\lambda_2\eta+\kappa_2)\}, \quad \forall (\tau,x)\in \mathcal{T}\times \mathbb{R}.
\label{interventionP2duop}
\end{align}

Recall that the Firm 1 and the Firm 2 profit functions are given by:
\begin{align}
&
\begin{aligned}
\Pi_1 (y;u,v)=\mathbb{E}^{[y]}  \Bigg[\int_{t_0}^{\tau_S}e^{-\epsilon r}[\alpha_1 X_1 (r)dr-\beta_1 X_2(r)]
&-\sum_{j\geq 1}c_1 (\tau_j,\xi_j ) \cdot 1_{\{\tau_j\leq\tau_S \}}\\& +\gamma_1e^{-\epsilon \tau_S} [X_1 (\tau_S )]^2 [X_2 (\tau_S )]^2 \Bigg], 
\end{aligned} 
\label{PIP1duop}
\\&
\begin{aligned}
\Pi_2 (y;u,v)=\mathbb{E}^{[y]}  \Bigg[\int_{t_0}^{\tau_S}e^{-\epsilon r}[\alpha_2 X_2 (r)-\beta_2 X_1(r)] dr
&-\sum_{m\geq 1}c_2 (\rho_m,\eta_m ) \cdot 1_{\{\rho_m\leq\tau_S \}} \\&+\gamma_2e^{-\epsilon \tau_S} [X_1 (\tau_S )]^2 [X_2 (\tau_S )]^2 \Bigg], 	
\end{aligned}
\label{PIP2duop}
\end{align}
where $x_i:=X_i(t_0)\in\mathbb{R}_{>0}$.

Given the setup of Theorem \ref{Verification_theorem_for_Non-Zero-Sum_Games_with_Impulse_Control}, at time $s\in\mathbb{R}_{>0}$ the Firm $i$ running cost $f_i$ is now given by: $f_i(Y(s))=e^{-\epsilon s}(\alpha_iX_i(s)-\beta_iX_j(s))$; $i,j\in\{1,2\}$, the Firm $i$ intervention costs are given by: $c_i(\tau,\xi)=\lambda_i\xi+\kappa_i$ for some intervention time $\tau\in\mathcal{T}$ and intervention $\xi\in\mathcal{Z}$ and the Firm $i$ terminal reward is given by: $G_i(Y(\tau_S))=\gamma_ie^{-\epsilon \tau_S}[X_i(\tau_S)]^2[X_j(\tau_S)]^2$.

W.l.o.g. we shall focus on the case for Firm 1, the arguments for Firm 2 being analogous. We can now apply the conditions of Theorem \ref{Verification_theorem_for_Non-Zero-Sum_Games_with_Impulse_Control} to show that the value function is a solution to the following Stefan problem:
\begin{align}
&\mathcal{L}\phi_i(y)+f_i=0,\hspace{29 mm} &\forall (x_1,x_2)\in D\equiv D_1\cap D_2,\; i\in\{1,2\}, \label{contequationduop}
\\
&\frac{\partial}{\partial z}\phi_1(x_1+z,x_2)=e^{-\epsilon t}\lambda_1, &\forall (x_1,x_2)\notin D\equiv D_1\cap D_2. \label{stefan2}
\end{align}
Indeed (\ref{contequationduop}) is immediately observed using (iv') of Theorem \ref{Verification_theorem_for_Non-Zero-Sum_Games_with_Impulse_Control}. This implies that:
\begin{align}\nonumber
0=&\alpha_1e^{-\epsilon t} x_1-\beta_1e^{-\epsilon t}x_2+ \frac{\partial \phi_1}{\partial t}(y) +\sum_{j=1}^{2}\mu_j\frac{\partial \phi_1}{\partial x_j}(y)+\frac{1}{2}\sum_{i,j=1}^2\sigma^2_{ij}\frac{\partial^2\phi_1(y)}{\partial x_i \partial x_j}
\\&+\int_\mathbb{R}\{\phi_1(s,x_1+\theta_{1j}z,x_2+\theta_{2j}z)-\phi_1(y)-\theta_{1j}z\frac{\partial \phi_1}{\partial x_1}(y)-\theta_{2j}z\frac{\partial \phi_1}{\partial x_2}(y)\}\nu_j(dz).    
\label{contequationgeneratorduop}\end{align}

We now try a candidate for the function, i.e. we specify the form:
\begin{equation}
\phi_1(y)=e^{-\epsilon t}\psi_1(x_1,x_2),
\label{candidateexp}
\end{equation}
for some $\epsilon>0$ and, as of yet, undetermined function $\psi\in\mathcal{C}^{1,2}$.

After plugging (\ref{candidateexp}) into (\ref{contequationgeneratorduop}), we find that:
\begin{align}
\nonumber
&\alpha_1 x_1-\beta_1x_2-\epsilon\psi_1(x_1,x_2)+\sum_{j=1}^2\mu_j\frac{\partial\psi_1(x_1,x_2)}{\partial x_j}+\frac{1}{2}\sum_{i,j=1}^2\sigma^2_{ij}\frac{\partial^2\psi_1(x_1,x_2)}{\partial x_i \partial x_j}
\\&\begin{aligned}
+\sum_{j=1}^2\int_\mathbb{R}\{\psi_1(x_1+\theta_{1j}z,x_2+\theta_{2j}z)-\psi_1(x_1,x_2)&-z\theta_{1j}\frac{\partial\psi_1(x_1,x_2)}{\partial x_1}&
\\&-z\theta_{2j}\frac{\partial\psi_1(x_1,x_2)}{\partial x_2}\}\nu_j(dz)=0.
\label{contequationgeneratorcandidateexpduop}
\end{aligned}
\end{align}

Let us now suppose that: 
\begin{align}
\psi_1(x_1,x_2)\equiv\hat{\psi}_1(x_1)+\hat{\psi}_1(x_2) \label{psi_decomp_ansatz_2_ch4}.    
\end{align} 
Hence using (\ref{contequationgeneratorcandidateexpduop}) we deduce that:
\begin{align}\nonumber
\alpha_1 x_1-\beta_1x_2&+\sum_{i=1}^2\Big\{-\epsilon\hat{\psi}_1(x_i)+\mu_i\frac{\partial\hat{\psi}_1(x_i)}{\partial x_i}+\frac{1}{2}\sum_{i=1}^2\sigma^2_{ii}\frac{\partial^2\hat{\psi}_1(x_i)}{\partial x_i^2 }\Big\}\\&+\sum_{i,j=1}^2\int_\mathbb{R}\{\hat{\psi}_1(x_i+\theta_{ij}z)-\hat{\psi}_1(x_i)-z\theta_{ij}\frac{\partial\hat{\psi}_1(x_i)}{\partial x_i}\}\nu_j(dz)=0. \label{contequationgeneratorcandidatelinearduop}    
\end{align}
After which we find that $\hat{\psi}_1$ is a solution to
\begin{align}
h(y)=A_1e^{r_1x_1}+A_2e^{r_2x_1}+\frac{\alpha_1}{\epsilon}x_1+B_1e^{r_1x_2}+B_2e^{r_2x_2}-\frac{\beta_1}{\epsilon}x_2+\frac{1}{2\epsilon^2}(\mu_1\alpha_1-\mu_2\beta_1),
\label{psihsolution}\end{align}
where $A_1,A_2,B_1,B_2\in\mathbb{R}$ are unknown constants and $r_1$ and $r_2$ are roots of the equation:
\begin{equation}
q(r_k):=\frac{1}{2}\sigma_{ii}^2r_k^2+\mu_i r_k-\epsilon+\int_{\mathbb{R} }\{e^{r_k\theta_{ij}z}-1-\theta_{ij}r_kz\}\nu_j(dz),\quad i,j,k\in\{1,2\}.\label{ch_4_r_roots}
\end{equation}

W.l.o.g. let us set $r_1<r_2$. Now since $\lim_{|r|\to\infty}q(r)=\infty$ $\mathbb{P}-$a.s., and $q(0)=-\epsilon<0$ and since $\forall r,z$ we have that: $\{e^{r_k\theta_{ij}z}-1-\theta_{ij}r_kz\}\nu_j(dz)>0$ we find that:
\begin{equation}
|r_1|>r_1,
\end{equation}
and
\begin{equation}
r_1<0<r_2.
\end{equation}
Our ansatz for the continuation region $D_1$ is that it takes the form:
\begin{equation}
D_1=\{x_1> x_1^{\star}|x_1,x_1^{\star}\in\mathbb{R}\}. \label{contregionansatzduop}
\end{equation}
We now derive (\ref{stefan2}) and in doing so we shall determine $x^{\star}_1$.
Now for all $x_1\leq x^{\star}_1$ we have that:
\begin{equation}
\psi_1(x_1,x_2)=\mathcal{M}_1\psi_1(x_1,x_2)=\sup_{z\in\mathcal{Z}}\{\psi_1(x_1+z,x_2)-(\kappa_1+\lambda_1z)\}.\label{interventioneqduop}
\end{equation}
We wish to determine the value $z$ that maximises (\ref{interventioneqduop}), hence let us now define the function $G$ by the following expression:
\begin{equation}
G(\xi)=\psi_1(x_1+\xi,x_2)-(\kappa_1+\lambda_1\xi), \label{functionG}
\end{equation}
$\forall  \xi\in\mathcal{Z}, x_1,x_2 \in \mathbb{R}$.

We now seek to evaluate the maxima of (\ref{functionG}), i.e. when:
\begin{equation}
G'(\xi)=0.
\end{equation}
We therefore see that the following expression holds :
\begin{equation}
\psi'_1(x_1+\xi,x_2)=\lambda_1, \qquad \forall x_1,x_2 \in \mathbb{R},\xi \in \mathcal{Z}. \label{psiderivativeeqn1}
\end{equation}

Let us now consider a unique point $\hat{x}_1 \in ]0,x_1^{\star}[$ then:
\begin{equation}
\psi'_1(\hat{x}_1,x_2)=\lambda_1. \label{psiderivativeeqn2}
\end{equation}
Hence, we have that
\begin{equation}
\hat{x}_1=x_1+\hat{\xi}(x_1) \text{ or } \hat{\xi}(x_1)=\hat{x}_1-x_1.
\end{equation}

We therefore deduce that for $x\in]0,x^{\star}_1[$, we have that:
\begin{equation}
\psi_1(x_1,x_2)=\psi_1(\hat{x}_1,x_2)-\kappa_1+\lambda_1(x_1-\hat{x}_1), \label{psieqn}
\end{equation}
or which by \eqref{psi_decomp_ansatz_2_ch4} may be equivalently expressed as:
\begin{equation}
\psi_1(x_1,x_2)=\hat{\psi}_1(\hat{x}_1)+\hat{\psi}_1(x_2)-\kappa_1+\lambda_1(x_1-\hat{x}_1). \label{psieqn_2}
\end{equation}
Using (\ref{psiderivativeeqn1}) - (\ref{psiderivativeeqn2}) and (\ref{psieqn_2}) and inserting (\ref{psihsolution}), we can construct the following system of equations:
\begin{align}
A_1r_1e^{r_1x_1^{\star}}&+A_2r_2e^{r_2x_1^{\star}}+\frac{\alpha_1}{\epsilon}=\lambda_1,\label{riccartiproof1}\\
A_1r_1e^{r_1\hat{x}_1}&+A_2r_2e^{r_2\hat{x}_1}+\frac{\alpha_1}{\epsilon}=\lambda_1,\label{riccartiproof2}\\
A_1(e^{r_1x_1^{\star}}-e^{r_1\hat{x}_1})&+A_2(e^{r_2x_1^{\star}}-e^{r_2\hat{x}_1})=-\kappa_1 +\Big(\lambda_1-\frac{\alpha_1}{\epsilon}\Big)(x^{\star}_1-\hat{x}_1).\label{riccartiproof3}
\end{align}

Repeating the above steps for $\phi_2$ leads to an analogous set of equations as  (\ref{riccartiproof1})-(\ref{riccartiproof3}) with $(A_1,A_2,x^{\star}_1,\hat{x}_1,\alpha_1,\lambda_1)$ replaced by $(B_1,B_2,x^{\star}_2,\hat{x}_2,\alpha_2,\lambda_2)$.

Now, since the system (\ref{PIP1duop}) - (\ref{PIP2duop}) is invariant under the transformations $\{1\leftrightarrow 2\}$ then we must have $A_1=B_1(:=C_1)$ and $A_2=B_2(:=C_2)$ (since (\ref{psihsolution}) must still be a solution to (iv) after the transformation $\{1\leftrightarrow 2\}$). Hence, we are left with a system of 6 unknowns $(C_1,C_2,x_1^{\star},\hat{x}_1,x_2^{\star},\hat{x}_2)$ and 6 equations. We can therefore uniquely determine the values $(C_1,C_2,x_1^{\star},\hat{x}_1,x_2^{\star},\hat{x}_2)$ --- this proves Theorem \ref{Theorem 4.3.}.
\end{refproof}
\begin{refproof}[Proof of Proposition \ref{firm_value_function_ch4}]
The proof follows immediately after combining equations (\ref{candidateexp}) (the value function ansatz) and (\ref{psihsolution}) (the general solution of the value function), together with the system of equations (\ref{riccartiproof1}) - (\ref{riccartiproof3}) (and exploiting symmetry in accordance with the above remarks). We therefore find that the value function for Firm $i$ is given by the following:
\begin{align*}
\phi_i(y)=e^{-\epsilon t}\Big\{C_1(e^{r_1x_1}+e^{r_1x_2})+C_2(e^{r_2x_1}+e^{r_2x_2})+\frac{\alpha_i}{\epsilon}x_1-\frac{\beta_i}{\epsilon}x_2+\frac{1}{2\epsilon^2}(\mu_i\alpha_i-\mu_j\beta_i)\Big\},&
 \\i\neq j, i,j\in\{1,2\},&
\end{align*}
where the constants $(C_1,C_2,x_1^{\star},\hat{x}_1,x_2^{\star},\hat{x}_2)$ are determined by the solutions to the following system of equations:
\begin{align}
{C}_1 r_1e^{r_1x_i^{\star}}&+{C}_2 r_2e^{r_2x_i^{\star}}+\frac{\alpha_i}{\epsilon}&&\hspace{0 mm}=\lambda_i,\label{riccartiproof1_no_jumps_end}\\
{C}_1 r_1e^{r_1\hat{x}_i}&+{C}_2 r_2e^{r_2\hat{x}_i}+\frac{\alpha_i}{\epsilon}&&\hspace{0 mm}=\lambda_i,\label{riccartiproof2_no_jumps_end}\\
{C}_1 (e^{r_1x_i^{\star}}-e^{r_1\hat{x}_i})&+{C}_2 (e^{r_2x_i^{\star}}-e^{r_2\hat{x}_i})&&\hspace{0 mm}=-\kappa_i +\Big(\lambda_i-\frac{\alpha_i}{\epsilon}\Big)(x^{\star}_i-\hat{x}_i)\label{riccartiproof_no_jumps_end}\\\nonumber
&&&\hfill\qquad\qquad\qquad\qquad\qquad\qquad\qquad\qquad i\in\{1,2\},
\end{align}
and $r_i$ are the roots of the equation 
\begin{equation}
q(r_k):=\frac{1}{2}\sigma_{ii}^2r_k^2+\mu_i r_k-\epsilon+\int_{\mathbb{R} }\{e^{r_k\theta_{ij}z}-1-\theta_{ij}r_kz\}\nu_j(dz),\label{ch_4_r_roots_end}
\end{equation}
for $i,j,k\in\{1,2\}$.
\end{refproof}
The transcendental nature of the system of equations (\ref{riccartiproof1_no_jumps_end}) - (\ref{riccartiproof_no_jumps_end}) means that the solutions to the constants $(C_1,C_2,x_1^{\star},\hat{x}_1,x_2^{\star},\hat{x}_2)$ cannot be expressed in closed form. Similarly, we cannot find closed solutions for the constants $r_1,r_2$ fir the integral equation (\ref{ch_4_r_roots_end}). Nonetheless, as the following result shows, if we restrict our attention to the case in which the market does not contain exogenous shocks, we can recover closed solutions to the constants $r_1,r_2$.  
% \subsection{The case without jumps \texorpdfstring{\textbf{($\mathbf{\theta_{ij}=0; i,j\in\{1,2\}}$)}}]}
\subsection{The case without jumps {\textbf{($\mathbf{\theta_{ij}=0; i,j\in\{1,2\}}$)}}]}

For the case in which the market contains no exogenous shocks using (\ref{ch_4_r_roots}), we readily observe that the constants $r_1,r_2$ can be solved analytically. In this case, each value function $\phi_i$ is given by
\begin{align*}
\phi_i(y)=e^{-\epsilon t}\left\{\tilde{C}_1(e^{r_1x_1}+e^{r_1x_2})+\tilde{C}_2(e^{r_2x_1}+e^{r_2x_2})+\frac{\alpha_i}{\epsilon}x_1-\frac{\beta_i}{\epsilon}x_2+\frac{1}{2\epsilon^2}(\mu_i\alpha_i-\mu_j\beta_i)\right\}&
\\i\neq j, i,j\in\{1,2\}&
\end{align*}
where the constants $(\tilde{C}_1,\tilde{C}_2,x_1^{\star},\hat{x}_1,x_2^{\star},\hat{x}_2)$ are determined by the solutions to the following system of equations:
\begin{align}
\tilde{C}_1 r_1e^{r_1x_i^{\star}}&+\tilde{C}_2 r_2e^{r_2x_i^{\star}}+\frac{\alpha_i}{\epsilon}&&\hspace{0 mm}=\lambda_i,\label{riccartiproof1_no_jumps}\\
\tilde{C}_1 r_1e^{r_1\hat{x}_i}&+\tilde{C}_2 r_2e^{r_2\hat{x}_i}+\frac{\alpha_i}{\epsilon}&&\hspace{0 mm}=\lambda_i,\label{riccartiproof2_no_jumps}\\
\tilde{C}_1 (e^{r_1x_i^{\star}}-e^{r_1\hat{x}_i})&+\tilde{C}_2 (e^{r_2x_i^{\star}}-e^{r_2\hat{x}_i})&&\hspace{0 mm}=-\kappa_i +\Big(\lambda_i-\frac{\alpha_i}{\epsilon}\Big)(x^{\star}_i-\hat{x}_i),\\
&&&\nonumber\hfill\qquad\qquad\qquad\qquad\qquad\qquad i\in\{1,2\},
\label{riccartiproof6_no_jumps}
\end{align}
where the constants $r_1$ and $r_2$ are given by:
\begin{equation}
r_1=-\frac{1}{\sigma_{ii}^2}\left(\mu_i+\sqrt{\mu_i^2+2\sigma_{ii}^2\epsilon}\right) ,\qquad r_2=\frac{1}{\sigma_{ii}^2}\left(\sqrt{\mu_i^2+2\sigma_{ii}^2\epsilon}-\mu_i\right).
\label{ch_4_r_roots_no_jumps}
\end{equation}

\section{Conclusion}

Using standard assumptions, we proved a verification theorem for a stochastic differential game with impulse controls, then generalising the results to cover a non-zero-sum payoff structure where the appropriate equilibrium concept is a Nash equilibrium. Having characterised the value for the stochastic differential game, we then applied the results to characterise optimal investment strategies for the dynamic advertising duopoly problem described in Section \ref{section_duopoly}.

% An interesting question for future research is investigating the above framework (in which the controllers use impulse controls to modify the state process) when either or both of the players only has access to partial information. Of particular interest is partial state information --- that is, a system in which the state process is adapted to some subset of the canonical filtration. A single player impulse controller version of a framework in which players have partial state information was studied in \cite{oksendal2008optimal}, here the controller's actions are subject to some execution delay so that there is some non-zero lag between the decision to apply an impulse intervention and the execution being carried out.
\section{Appendix}
A.1.1. Lipschitz Continuity

We assume there exist real-valued constants $c_{\mu},c_{\sigma}>0$ and $c_{\gamma} (\cdot)\in L^1\cap L^2 ( \mathbb{R}^l,\nu)$ s.th. $\forall  s\in [0,T], \forall  x,y\in S$ and $\forall  z\in \mathbb{R}^l$ we have:
\begin{align*}
|\mu(s,x)-\mu(s,y)|&\leq c_{\mu} |x-y|\\
|\sigma(s,x)-\sigma(s,y)|&\leq c_{\sigma} |x-y|
\\\int_{|z|\geq 1} |\gamma(x,z)-\gamma(y,z)| &\leq c_{\gamma} (z)|x-y|.
\end{align*}
A.1.2. Lipschitz Continuity

We also assume the Lipschitzianity of the functions $f$ and $g$ that is, we assume the existence of real-valued constants $c_f,c_g>0$ s.th. for $R\in \{f,g\}$:
\begin{equation*}
|R(s,x)+R(s,y)|\leq c_R |x-y|,\qquad \forall  s\in [0,T], \forall  (x,y) \in S.
\end{equation*}
A.2. Growth Conditions

We assume the existence of a real-valued constants $d_{\mu},d_{\sigma}>0$ and $d_{\gamma} (\cdot)\in L^1\cap L^2 ( \mathbb{R}^l,\nu), \rho\in [0,1)$ s.th. $\forall  (s,x)\in [0,T]\times S$ and $\forall  z\in \mathbb{R}^l$  we have:
\begin{align*}
|\mu(s,x)|\leq d_{\mu} (|1+|x|^\rho |)
\\
|\sigma(s,x)|\leq d_{\sigma} (|1+|x|^\rho |)
\\
\int_{|z|\geq 1} |\gamma(x,z)| \leq d_{\gamma} (|1+|x|^\rho |).
\end{align*}

Assumptions A.1.1 and A.2 ensure the existence and uniqueness of a solution to (\ref{uncontrolledstateprocess}) (c.f. \cite{cosso2013stochastic}). 
% Assumption A.3 (i) (subadditivity) is required in the proof of the uniqueness of the value function. Assumption A.3 (ii) (the player cost function is a decreasing function in time) and may be interpreted as a discounting effect on the cost of interventions.
Assumption A.1.2 is required to prove the regularity of the value function (see for example \cite{bensoussan1982controle} and for the single-player case, see for example \cite{chen2013impulse}). 
% Assumption A.3 (ii) was introduced (for the two-player case) in \cite{stroock2007multidimensional} though is common in the treatment of single-player case problems (e.g. \cite{chen2013impulse,tang1993finite}). Assumption A.4 is integral to the definition of the impulse control problem.\\
\section{Technical Conditions: (T1) - (T4)}
\renewcommand{\theenumi}{\roman{enumi}}
 \begin{enumerate}[leftmargin= 6 mm]

\item[(T1)] Assume that $\mathbb{E}[\int_{0}^{T} 1_{\partial D} (X^{\cdot,u} (s))ds]=0$ for all $x \in S,\; \forall u \in \mathcal{U}$ where $D\equiv D_1\cup D_2$.

\item[(T2)] $\partial D$ is a Lispchitz surface, that is to say that $\partial D$ is locally the graph of a Lipschitz continuous function: $\phi \in \mathcal{C}^2 (S\backslash \partial D) $ with locally bounded derivatives.

\item[(T3)] The sets $\{\phi^- (X^{\cdot,u} (\tau_m ));$ $ \tau_m \in \mathcal{T},\forall  m \in \mathbb{N}\}$ and $\{\phi^- (X^{\cdot,u} (\rho_j));\rho_j \in \mathcal{T},\forall  j \in \mathbb{N}\}$ are uniformly integrable $\forall  x \in S$ and $\forall u \in \mathcal{U}$.

	\item[(T4)] $\mathbb{E}[|\phi(X^{\cdot,u} (\tau ))|+|\phi(X^{\cdot,u} (\rho))|+\int_{0}^{T} |\mathcal{L}\phi(X^{\cdot,u} (s))|ds]< \infty,$ $ \\   \text{for all intervention times } \tau_,\rho \in \mathcal{T}$ and $\forall u \in \mathcal{U}.$
\end{enumerate}
Assumptions T3 and T4 hold if for example $\phi$ satisfies a polynomial growth condition and guarantee that $\phi(X(\tau)$ is both well-defined and finite --- in particular the uniform integrability assumption (T3) implies the existence of a finite constant $c>0$ s.th. $\mathbb{E}[|\phi(X(\tau))|1_{\{\tau<\infty\}}]\leq c$ for all $\tau\in\mathcal{T}$.
\begin{refproof}[Proof of Lemma \ref{Lemma 3.1.6.}]

Firstly, let us set $\Theta_i(s)\equiv \lambda_i, i\in \{1,2\}$ and suppose that $\xi(s)\equiv \nu_1^+(s)-\nu_1^- (s)$ and $\eta(s)\equiv\nu_2^+(s)-\nu_2^-(s)$ for player I and player II controls (resp.) where $s\in[0,T]$ and $\nu_i^+$ and $\nu_i^-$, $i \in \{1,2\}$ are given by the using expressions:
\begin{align}
\nu_1^+ (s)&=\frac{1}{2}\Bigg[\sum_{j\geq 1}(\xi_{j}\cdot 1_{\{\xi_j>0\}}+\lambda_1^{-1}\kappa_1\cdot 1_{\{\tau_j\leq s\}}-\lambda_1^{-1}\kappa_1)  \cdot 1_{\{\tau_j\leq s\}} \Bigg], \label{p1uppersingdefn}\\
\nu_1^-(s)&=- \frac{1}{2}\Bigg[\sum_{j\geq 1}(\xi_j\cdot 1_{\{\xi_j<0\}}+\lambda_1^{-1}\kappa_1\cdot 1_{\{\tau_j\leq s\}}-\lambda_1^{-1}\kappa_1 ) \cdot 1_{\{\tau_j\leq s\}}\Bigg] 	  \label{p1lowersingdefn}
\end{align}
and similarly for the player II control:
\begin{align}
\nu_2^+ (s)&=\frac{1}{2}\Bigg[\sum_{m\geq 1}( \eta_m\cdot 1_{\{\eta_m>0\}}+\lambda_2^{-1}\kappa_2\cdot 1_{\{\rho_m\leq s\}}-\lambda_2^{-1}\kappa_2 ) \cdot 1_{\{\rho_m\leq s\}}\Bigg] \label{p2uppersingdefn}\\
\nu_2^- (s)&=- \frac{1}{2}\Bigg[\sum_{m\geq 1}(\eta_m\cdot 1_{\{\eta_m<0\}}+\lambda_2^{-1}\kappa_2\cdot 1_{\{\rho_m\leq s\} }-\lambda_2^{-1}\kappa_2 ) \cdot 1_{\{\rho_m\leq s\} } \Bigg]
\label{p2lowersingdefn}
\end{align}

We do the proof for the player II impulse controls,  the proof for the player I part is analogous. 

Using (\ref{p1uppersingdefn}) - (\ref{p1lowersingdefn}) we readily deduce that:
\begin{align*}
d\eta(s)&=d\nu_2^+(s)-d\nu_2^- (s)
\\&
=\frac{1}{2}\sum_{m\geq 1}( \eta_m\cdot 1_{\{\eta_m>0\}}+2\lambda_2^{-1}\kappa_2\cdot 1_{\{\rho_m\leq s\}}-\lambda_2^{-1}\kappa_2 ) \cdot \delta_{\rho_m}(s)
\\&\qquad
+ \frac{1}{2}\sum_{m\geq 1}( \eta_m\cdot 1_{\{\eta_m<0\}}+2\lambda_2^{-1}\kappa_2\cdot 1_{\{\rho_m\leq s\}}-\lambda_2^{-1}\kappa_2 ) \cdot \delta_{\rho_m}(s)
\\&
=\sum_{m\geq 1}( \eta_m+(2\cdot1_{\{\rho_m\leq s\}}-1)\lambda_2^{-1}\kappa_2)\cdot \delta_{\rho_m}(s)
\end{align*}

Using the properties of the Dirac-delta function and by Fubini's theorem we find that:
\begin{align*}
\int_{t_0}^{\tau_S}\Theta_2(s) d\eta(s)&=\sum_{m\geq 1}\int_{t_0}^{\tau_S}(\lambda_2\eta_m+(2\cdot1_{\{\rho_m\leq s\}}-1)\kappa_2)\cdot \delta_{\rho_m}(s)
\\&=\sum_{m\geq 1}\int_{t_0}^{\tau_S}(\lambda_2\eta_m+(2\cdot1_{s\in [\rho_m,\infty[}-1)\kappa_2)\cdot \delta_{\rho_m}(s)
\\&=\sum_{m= 1}^{\mu_{[t_0,\tau_S]}(v)}(\lambda_2\eta_m+\kappa_2)
\\&=\sum_{m \geq 1}\chi (\rho_m,\eta_m )\cdot 1_{\{\rho_m\leq\tau_S\}}. \end{align*}

Lastly, we compute $\eta(s)$, indeed we observe that:
\begin{align*}
\eta(s)&=\nu_2^+(s)+\nu_2^- (s)	
\\&
=\sum_{m\geq 1}(( \eta_m+\lambda_2^{-1}\kappa_2\cdot 1_{\{\rho_m\leq s\}})-\lambda_2^{-1}\kappa_2 ) \cdot 1_{\{\rho_m\leq s\}}\cdot 1_{\{\eta_m>0\}}
\\&\qquad+ \sum_{m\geq 1}(( \eta_m+\lambda_2^{-1}\kappa_2\cdot 1_{\{\rho_m\leq s\}})-\lambda_2^{-1}\kappa_2 ) \cdot1_{\{\rho_m\leq s\}}\cdot 1_{\{\eta_m<0\}}.
\end{align*}
Now, since $1_{\{\rho_m\leq s\}}\cdot 1_{\{\rho_m\leq s\}}= 1_{\{\rho_m\leq s\}}$ we find that:
\begin{align*} 
\eta(s)=\sum_{m\geq 1}(\eta_m+\lambda_2^{-1}\kappa_2-\lambda_2^{-1}\kappa_2)\cdot 1_{\{\rho_m\leq s\}}
=\sum_{m\geq 1}\eta_m\cdot 1_{\{\rho_m\leq s\}}
=\sum_{m= 1}^{\mu_{[t_0,s]}(v)}\eta_m.\end{align*}

Hence, after repeating the exercise for the player I controls (using that $\xi(s):=\nu_1^+(s)-\nu_1^- (s)$) and setting $\nu(s)=\nu^+(s)-\nu^-(s)$ where $\nu^+(s)\equiv \nu^+_1(s) +\nu^+_2(s)$, $\nu^-\equiv \nu^-_1 +\nu^-_2$  we recover the impulse control game.
\end{refproof}
\begin{refproof}[Proof of Theorem \ref{Verification_theorem_for_Non-Zero-Sum_Games_with_Impulse_Control}]

We prove the theorem for player I with the proof for player II being analogous.

As in the proof of Theorem \ref{Verification_theorem_for_Zero-Sum_Games with Impulse Control}, we begin by adopting the following notation:
\begin{align}
    &Y^{y_0,\cdot}(s)\equiv (s,X^{t_0,x_0,\cdot}(t_0+s)), \quad y_0\equiv (t_0,x_0), \; \forall s\in [0,T-t_0], \\& \hat{Y}^{y_0,\cdot}(\tau)=Y^{y_0,\cdot}(\tau^{-} )+\Delta_N Y^{y_0,\cdot}(\tau), \quad \tau\in\mathcal{T},
\end{align}
where $\Delta_N Y(\tau) $ denotes a jump at time $\tau$ due to $\tilde{N}$.

Correspondingly, we adopt the following impulse response function $\hat{\Gamma}: \mathcal{T}\times S\times \mathcal{Z}\to  \mathcal{T}\times S$ acting on  $y'\equiv (\tau,x')\in \mathcal{T}\times S$ where $x'\equiv X^{t_0,x_0,\cdot}(t_0+\tau^-)$ is given by: 
\begin{align}
\hat{\Gamma}(y',\zeta)\equiv (\tau,\Gamma (x',\zeta))=(\tau,X^{t_0,x_0,\cdot} (\tau)),\quad \forall \xi\in\mathcal{Z},\; \forall\tau\in\mathcal{T} .
\end{align}

Let us also now fix the player II control $\hat{v}\in \mathcal{V}$; we firstly appeal to Dynkin's formula for jump-diffusion processes hence, we have the following:
\begin{align*}
\mathbb{E}[\phi_1(Y^{y_0,u,\hat{v}} (\tau_j ))]-\mathbb{E}[\phi_1(\hat{Y}^{y_0,u,\hat{v}} (\tau_{j+1}-))]=-\mathbb{E}\left[\int_{\tau_j}^{\tau_{j+1}}\frac{\partial \phi_1}{\partial s}+\mathcal{L}\phi(Y^{y_0,u,\hat{v}} (s))ds\right].  \label{dynkin1NEproof} \end{align*}
Summing from $j=0$ to $j=k$  implies that:
\begin{align}
\phi_1 (y_0)+\sum_{j=1}^k&\mathbb{E}\left[\phi_1 (\hat{Y}^{y_0,u,\hat{v}} (\tau_j ))-\phi_1 (\hat{Y}^{y_0,u,\hat{v}} (\tau_j^- ))\right] -\mathbb{E}\left[\phi_1 (\hat{Y}^{y_0,u,\hat{v}} (\tau_{k+1}^- ))\right] 	
\\&=-\mathbb{E}\left[\int_{t_0}^{\tau_{k+1}}\frac{\partial \phi_1}{\partial s}+\mathcal{L}\phi_1(Y^{y_0,u,\hat{v}} (s))ds\right].   	\label{dynkin1NEproofSUM} \end{align}
Now by similar reasoning as in the zero-sum case (c.f. (\ref{interventionineq2theorem6.3})), we have that:
\begin{align}\nonumber
&\mathcal{M}_1 \phi_1 (\hat{Y}^{y_0,u,\hat{v}} (\tau_j ^-))-\phi_1 (\hat{Y}^{y_0,u,\hat{v}} (\tau_j^- ))+c_1 (\tau_j,\xi_j )\cdot1_{\{\tau_j\leq \tau_S\}}
\\&\geq \phi_1 (Y^{y_0,u,\hat{v}} (\tau_j ))-\phi_1 (\hat{Y}^{y_0,u,\hat{v}} (\tau_j^- )).\label{interventionineqNEproof}\end{align} 
Inserting (\ref{interventionineqNEproof}) into (\ref{dynkin1NEproofSUM}) implies that:
\begin{align}
&\begin{aligned}\phi_1 (y_0)+\sum_{j=1}^k\mathbb{E}[\mathcal{M}_1 \phi_1 (\hat{Y}^{y_0,u,\hat{v}} (\tau_j ^-))-\phi_1 (\hat{Y}^{y_0,u,\hat{v}} (\tau_j^- ))&+c_1 (\tau_j,\xi_j )\cdot1_{\{\tau_j\leq \tau_S\}}
\\&-\mathbb{E}[\phi_1 (\hat{Y}^{y_0,u,\hat{v}} (\tau_{k+1}^- ))] 	\nonumber
\end{aligned}\\&\geq-\mathbb{E}\left[\int_{t_0}^{\tau_{k+1}}\frac{\partial \phi_1}{\partial s}+\mathcal{L}\phi_1(Y^{y_0,u,\hat{v}} (s))ds\right].    \label{continuityrunningpostsumNE}\end{align}
Additionally, by (ii') we have that:
\begin{align*}
&\qquad\qquad\qquad\frac{\partial \phi}{\partial s}+\mathcal{L}\phi_1 (Y^{y_0,u_{[t_0,s]} ,\hat{v}_{[t_0,s]}} (s)) 	
\\&\geq \frac{\partial \phi}{\partial s}+\mathcal{L}\phi_1 (Y^{\hat{u}_{[t_0,s]} ,\hat{v}_{[t_0,s]}} (s))+f_1 (Y^{y_0,\hat{u}_{[t_0,s]} ,\hat{v}_{[t_0,s]}} (s))-f_1 (Y^{y_0,u_{[t_0,s]} ,\hat{v}_{[t_0,s]}} (s))) 	
\\&\geq -f_1 (Y^{y_0,u_{[t_0,s]} ,\hat{v}_{[t_0,s]}} (s)). 
\end{align*}
Or
\begin{align}
-\left(\frac{\partial \phi}{\partial s}+\mathcal{L}\phi_1 (Y^{y_0,u_{[t_0,s]} ,\hat{v}_{[t_0,s]}} (s))\right)\leq f_1 (Y^{y_0,u_{[t_0,s]} ,\hat{v}_{[t_0,s]}} (s)). 	 \label{continuityrunningconditionNE}\end{align}
Hence, inserting (\ref{continuityrunningconditionNE}) into (\ref{continuityrunningpostsumNE}) yields:
\begin{align*}
&\begin{aligned}
\phi_1 (y_0)+ \sum_{j=1}^{k} \mathbb{E}[\mathcal{M}_1  \phi_1  (\hat{Y}^{y_0  ,u, \hat{v} } (\tau_j^-  ))- \phi_1  (\hat{Y}^{y_0  ,u, \hat{v} } (\tau_j ^-  ))&+c_1(\tau_j  , \xi_j  )\cdot 1_{\{\tau_j  \leq  \tau_S  \}} ] 
\\&- \mathbb{E}[\phi_1  (\hat{Y}^{y_0  ,u, \hat{v} } (\tau_{k+1}^-  ))]
\end{aligned}
\\&\geq \mathbb{E}\left[\int_{t_0}^{\tau_{k+1}}f_1 (Y^{y_0,u_{[t_0,s]} ,\hat{v}_{[t_0,s]}} (s))ds\right].
\end{align*}
Now, as in the proof for the zero-sum case, we have, using (ii) that\\$\lim_{s\to \tau_S^-}[\mathcal{M}_i \phi_i (Y^{y_0,u,v} (s))-\phi_i (\hat{Y}^{y_0,u,v} (s))]=0,  i\in \{1,2\}$ and\\$\lim_{s\to \tau_S^-}\phi(Y^{\cdot,u,v} (s))=G(Y^{\cdot,u,v} (\tau_S ))] \hspace{1 mm}\forall  u \in \mathcal{U},\; \forall v \in \mathcal{V},\;  \mathbb{P}-$a.s.. Hence, after taking the limit $k\to \infty$, we recognise:
\begin{align}
\phi_1 (y_0)\geq \mathbb{E}\Bigg[\int_{t_0}^{\tau_S}f_1 (Y^{y_0,u_{[t_0,s]} ,\hat{v}_{[t_0,s]}} (s))ds&-\sum_{j\geq1}c_1(\tau_j,\xi_j ) \cdot 1_{\{\tau_j\leq \tau_S \}}
\\&+G_1(Y^{y_0,u,v} (\tau_S ))\cdot 1_{\{\tau_S<\infty\}}  \Bigg],
\end{align}
or
\begin{equation}
\phi_1 (y)\geq J_1^{(u,\hat{v} )} [y], \quad \forall y\in [0,T]\times S.
\end{equation}
 Since this holds for all $u \in \mathcal{U}$ we have
 \begin{equation}
\phi_1 (y)\geq \sup_{u \in \mathcal{U}}J_1^{(u,\hat{v} )} [y], \quad \forall y\in [0,T]\times S.	\end{equation}
Now, applying the above arguments and fixing the pair of controls $(\hat{u},\hat{v} )\in \mathcal{U}\times \mathcal{V}$ yields the following equality:
\begin{equation}
\phi_1 (y)=\sup_{u \in \mathcal{U}}J_1^{(u,\hat{v} )} [y]=J_1^{\hat{u},\hat{v} } [y],\; {\forall  y\in [0,T]\times S.}
\end{equation}
After using an analogous argument for the player II policy $v\in \mathcal{V}$ (fixing the player I control as $\hat{u}$), we deduce that:
\begin{equation}
\phi_2 (y)\geq \sup_{u \in \mathcal{U}}J_2^{(u,\hat{v} )}  [y],\quad \forall y\in [0,T]\times S,
\end{equation}
after which we observe the following statements:
\begin{align}
\phi_2 (y)=\sup_{v \in \mathcal{V}}J_2^{(\hat{u},v)} [y]=J_2^{\hat{u},\hat{v} } [y],\\
\phi_1 (y)=\sup_{u \in \mathcal{U}}J_1^{(u,\hat{v} )} [y]=J_1^{\hat{u},\hat{v} } [y], 	
\end{align} 
from which we deduce that $(\hat{u},\hat{v})$ is a Nash equilibrium and hence the thesis is proven.
\end{refproof}
\end{example}

\printbibliography[
heading=bibintoc,
title={References}
]

\end{document}